# CONSISTENCY OF THE JACKKNIFE-AFTER-BOOTSTRAP VARIANCE ESTIMATOR FOR THE BOOTSTRAP QUANTILES OF A STUDENTIZED STATISTIC[1]

By S. N. Lahiri

*Iowa State University*

Efron [*J. Roy. Statist. Soc. Ser. B* **54** (1992) 83–111] proposed a computationally efficient method, called the jackknife-after-bootstrap, for estimating the variance of a bootstrap estimator for independent data. For dependent data, a version of the jackknife-after-bootstrap method has been recently proposed by Lahiri [*Econometric Theory* **18** (2002) 79–98]. In this paper it is shown that the jackknife-after-bootstrap estimators of the variance of a bootstrap quantile are consistent for both dependent and independent data. Results from a simulation study are also presented.

**1. Introduction.** Bootstrap quantiles are frequently used for the purpose of constructing bootstrap-based confidence intervals (CIs). Like any other estimation procedure, the accuracy of the quantile estimators produced by the bootstrap method needs to be assessed. In his seminal paper, Efron [11] formulated the jackknife-after-bootstrap (JAB) method for estimating the standard errors of bootstrap quantities in the independent case. Recently, Lahiri [27] proposed a version of the JAB method for assessing the accuracy of *block* bootstrap estimators in the dependent case. In this paper we establish theoretical validity of the JAB method for assessing the variability of bootstrap quantiles for both dependent and independent data. To describe the main results briefly, suppose that $\theta$ is a population parameter and $\hat{\theta}_n$ is an estimator of $\theta$ that can be represented as a smooth function of sample means. Let $\hat{G}_n$ denote the block bootstrap estimator of the distribution function of a *Studentized* version of $\hat{\theta}_n$. See Section 2 for more details. Also, for a given $\alpha \in (0,1)$, let $\hat{G}_n^{-1}(\alpha)$ denote the $\alpha$-quantile of $\hat{G}_n$. The main results of the paper show that, under some regularity conditions, and for

Received August 2001; revised October 2004.
[1]Supported in part by NSF Grants DMS-95-05124 and DMS-00-72571.

both dependent as well as independent data,

$$\widehat{\mathrm{Var}}_{\mathrm{JAB}}(\hat{\varphi}_n)/\mathrm{Var}(\hat{\varphi}_n) \xrightarrow{p} 1 \quad \text{as } n \to \infty \tag{1.1}$$

for $\hat{\varphi}_n = \hat{G}_n(x)$ and $\hat{\varphi}_n = \hat{G}_n^{-1}(\alpha)$, where $\widehat{\mathrm{Var}}_{\mathrm{JAB}}(\hat{\varphi}_n)$ is the JAB-estimator of the variance of $\hat{\varphi}_n$. Thus, the JAB variance estimators of the bootstrap distribution function and the bootstrap quantiles are consistent in the Studentized case.

For independent data, $\widehat{\mathrm{Var}}_{\mathrm{JAB}}(\hat{\varphi}_n)$ is defined using the standard jackknife variance estimator formula (cf. [11]). In the dependent case, however, the definition of $\widehat{\mathrm{Var}}_{\mathrm{JAB}}(\hat{\varphi}_n)$ is slightly different from the usual block-jackknife variance estimator of $\hat{\phi}_n$. Here, instead of deleting blocks of observations directly, one deletes "blocks of bootstrap blocks" in the definition of the JAB variance estimator. This modification distinguishes the JAB in the dependent case from its i.i.d. counterpart and is crucial for computational efficiency (cf. [27]). Nonetheless, we show that consistency of the JAB variance estimator continues to hold in the dependent case.

In the literature the JAB method has been employed in many important problems (cf. [8, 13, 21]), but without any theoretical justification. The main results of the paper provide theoretical justification for using the JAB method for estimating the variance of an important class of nonsmooth functionals of the bootstrap distribution for both dependent and independent data. Together with the results of Lahiri [27], they show that the JAB method can be effectively used not only to assess the variability of simple functionals, like the bootstrap bias and variance estimators, but also to derive consistent estimates of standard errors for more complicated functionals, like bootstrap quantiles. These results are particularly important because of the extensive use of bootstrap quantiles in various bootstrap CIs, whose coverage accuracy critically depends on the accuracy of the quantile estimators, especially in the dependent case. Another notable application of the JAB method is recently given by Lahiri, Furukawa and Lee [29] who make use of the JAB variance estimators to construct a data-based method for selecting the optimal block size for block bootstrap estimation of various functionals, including block bootstrap quantiles.

The main steps involved in the proofs of the main results are the derivations of certain *uniform* Edgeworth and Cornish–Fisher expansions for the jackknife pseudo-values of the (block) bootstrap distribution function and quantile estimators, respectively. The technical arguments developed to prove these uniform expansions for the Studentized statistics in the dependent case may be of some independent interest. We also report the results from a small simulation study on finite sample properties of the JAB variance estimators in the dependent case.

Properties of the jackknife and the bootstrap methods under independence have been studied by several authors; see Efron [9, 10], Miller [31],

Parr and Schucany [34], Shao and Tu [38], Shao and Wu [39] and the references therein. In his 1992 seminal paper, Efron formulated the JAB method for independent data and established its computational efficacy. The JAB method for dependent data was formulated by Lahiri [27]. Block bootstrap and block jackknife methods for dependent data have been put forward by several authors, notably by Hall [18], Carlstein [7], Künsch [22], Liu and Singh [30], Politis and Romano [35, 37] and Paparoditis and Politis [33]. Properties of block-bootstrap methods have been studied by Lahiri [23, 25, 26], Bühlmann [5], Naik-Nimbalkar and Rajarshi [32], Hall, Horowitz and Jing [20] and Götze and Künsch [17], among others.

The rest of the paper is organized as follows. In Section 2 we describe the basic framework of the paper and the JAB method for dependent data. Assumptions and the results on consistency of the JAB variance estimators are given in Section 3. Results from a simulation study are reported in Section 4. Proofs of the main results are given in Section 5.

**2. Preliminaries.** For brevity, we concentrate our discussion on the dependent case and, as appropriate, point to the main differences for the independent case. In Section 2.1 we describe the basic framework of this paper and in Section 2.2 the JAB method.

2.1. *The basic framework.* Let $X_1, X_2, \ldots$ be a sequence of stationary random vectors (r.v.s) with values in $\mathbb{R}^d$. For proving the results, we shall work under the "smooth function model" of Bhattacharya and Ghosh [1] (cf. [19]). Let $\theta$ be a parameter of interest and let $\hat{\theta}_n$ be an estimator of $\theta$. Let $\bar{X}_n = n^{-1} \sum_{i=1}^{n} X_i$ and $\mu = EX_1$. Then, the "smooth function model" prescribes that there exists a smooth function $H : \mathbb{R}^d \to \mathbb{R}$ such that $\hat{\theta}_n = H(\bar{X}_n)$ and $\theta \equiv H(\mu)$. Considering suitable transformations of the observations, one can treat a wide class of common estimators, including the generalized $M$-estimators of Bustos [6], under this model. See, for example, [19] and [28] for the scope of the "smooth function model" in the independent and the dependent cases, respectively.

Note that, under the smoothness assumption on the function $H$, by the delta method,

$$n^{1/2}(\hat{\theta}_n - \theta) \xrightarrow{d} N(0, \tau_\infty^2),$$

where $\xrightarrow{d}$ denotes convergence in distribution and $\tau_\infty^2 \equiv \lim_{n \to \infty} n \operatorname{Var}(h(\mu)' \bar{X}_n)$ denotes the asymptotic variance of $\hat{\theta}_n$. Here $h$ denotes the $d \times 1$ vector of first-order partial derivatives of $H$. Let $\hat{\Gamma}_n(k) = n^{-1} \sum_{i=1}^{n-k} (X_i - \bar{X}_n)(X_{i+k} - \bar{X}_n)'$ and let $\hat{\Sigma}_n = \hat{\Gamma}_n(0) + \sum_{k=1}^{\ell_1} \{\hat{\Gamma}_n(k) + \hat{\Gamma}_n(k)'\}$ for some integer $\ell_1$. Then we may use $\hat{\tau}_n^2 \equiv |h(\bar{X}_n)'\hat{\Sigma}_n h(\bar{X}_n)|$ as an estimator of $\tau_\infty^2$ and define the Studentized version of $\hat{\theta}_n$ as

(2.1) $$T_n \equiv n^{1/2}(\hat{\theta}_n - \theta)/(\hat{\tau}_n + n^{-1}).$$

In the dependent case $\ell_1$ goes to infinity suitably with the sample size $n$, while in the independent case we set $\ell_1 \equiv 0$ for all $n$. We add the factor $n^{-1}$ in the denominator of (2.1) to make $T_n$ well defined in the case when $\hat{\tau}_n = 0$.

The sampling distribution and the quantiles of the Studentized statistic $T_n$ may be estimated by the bootstrap methods. For the dependent case, here we restrict attention to the moving block-bootstrap (MBB) method of Künsch [22] and Liu and Singh [30]. (See remarks at the end of Section 3 for other block bootstrap methods.) Given $\mathcal{X}_n = (X_1, \ldots, X_n)$, define the blocks $\mathcal{B}_j = (X_j, \ldots, X_{j+\ell-1})$, $j = 1, \ldots, N$, where $N = n - \ell + 1$ and $1 \leq \ell \leq n$ denotes the block size. Let $b = \lfloor n/\ell \rfloor$ be the largest integer not exceeding $n/\ell$. Then for the MBB method one selects a random sample $\mathcal{B}_1^*, \ldots, \mathcal{B}_b^*$ from the collection $\{\mathcal{B}_1, \ldots, \mathcal{B}_N\}$. Arranging the components of $\mathcal{B}_1^*, \ldots, \mathcal{B}_b^*$ into a sequence yields $n_1 \equiv b \cdot \ell$ bootstrap samples $\mathcal{X}_n^* \equiv \{X_1^*, \ldots, X_{n_1}^*\}$. For $\ell = o(n)$ (which we assume in this paper), we get $n_1/n \to 1$ as $n \to \infty$. Let $\bar{X}_n^* = n_1^{-1} \sum_{i=1}^{n_1} X_i^*$ denote the bootstrap sample mean and let $\hat{\mu}_n = E_* \bar{X}_n^*$, where $E_*$ denotes the conditional expectation given $\{X_n\}_{n \geq 1}$. Then the bootstrap version of $(\hat{\theta}_n - \theta)$ is given by $(\theta_n^* - \tilde{\theta}_n)$ with $\theta_n^* = H(\bar{X}_n^*)$ and $\tilde{\theta}_n = H(\hat{\mu}_n)$. It is known (cf. [23]) that the "naive" approach of centering $\theta_n^*$ at $\hat{\theta}_n$ leads to a worse performance of the resulting bootstrap approximation compared to centering $\theta_n^*$ at $\tilde{\theta}_n$.

Next, to define the MBB version of $T_n$ [cf. (2.1)], note that, conditional on $\mathcal{X}_n$, the resampled block sums $S_i^* \equiv \sum_{j=(i-1)\ell+1}^{i\ell} X_j^*$ are i.i.d. Hence, following Götze and Künsch [17], we may use the sample variance of these block sums (with a suitable scaling) to define the Studentized version of $\theta_n^*$. A "natural" analog of $\hat{\tau}_n^2$ for the bootstrap case is thus given by

$$\tau_n^{*2} \equiv \frac{1}{\ell b} \sum_{i=1}^{b} \{h(\bar{X}_n^*)'(S_i^* - \ell \bar{X}_n^*)\}^2.$$

Since the MBB resamples blocks *with replacement* from the observed blocks $\mathcal{B}_1, \ldots, \mathcal{B}_N$, it is possible that $\tau_n^{*2} = 0$ for certain values of $X_n^*, \ldots, X_{n_1}^*$ (e.g., when all the resampled blocks $\mathcal{B}_1^*, \ldots, \mathcal{B}_b^*$ are equal to $\mathcal{B}_1$, say). To avoid division by zero, we define the bootstrap version of $T_n$ as

(2.2) $$T_n^* \equiv n_1^{1/2}(\theta_n^* - \tilde{\theta}_n)/(\tau_n^* + n^{-1}).$$

The inclusion of the term $n^{-1}$ makes $T_n^*$ well defined when $\tau_n^* = 0$. The effect of this modification on the accuracy of the bootstrap distribution function estimator is negligible up to the second order.

Next note that setting the block size $\ell \equiv 1$ for all $n$, the MBB reduces to the standard bootstrap method of Efron [9] and $T_n^*$ of (2.2) also gives the bootstrap version of $T_n$ in the independent case. For both the dependent and the independent cases, the bootstrap estimators of the unknown distribution function and the quantiles of $T_n$ are, respectively, given by

$$(2.3) \qquad \hat{\varphi}_{1n} \equiv P_*(T_n^* \leq x), \qquad x \in \mathbb{R}$$

and

$$(2.4) \qquad \hat{\varphi}_{2n} \equiv \inf\{x : P_*(T_n^* \leq x) \geq \alpha\}, \qquad \alpha \in (0,1),$$

where $P_*$ denotes conditional probability given $\{X_n\}_{n \geq 1}$. Variances of these bootstrap estimators can be estimated by the JAB method, which we describe next.

2.2. *Jackknife-after-bootstrap.* The JAB method for dependent data applies a version of the block jackknife method to a block-bootstrap estimator. To describe the (block) jackknife method, let $\hat{\gamma}_n = g_n(\mathcal{X}_n)$ be an estimator of a parameter $\gamma$. For $i = 1, \ldots, N$, define the $i$th jackknife block-deleted point value $\hat{\gamma}_n^{(i)}$ by $\hat{\gamma}_n^{(i)} = g_{n-\ell}(\mathcal{X}_{n,i})$, where $\mathcal{X}_{n,i} \equiv \mathcal{X}_n \setminus \mathcal{B}_i$ denotes the set of data values after deleting the $i$th block. The *block jackknife* estimator of the variance of $\hat{\gamma}_n$, *based on blocks of size* $\ell$, is given by (cf. [22, 30])

$$(2.5) \qquad \widehat{\mathrm{Var}}_{\mathrm{BJ}}(\hat{\gamma}_n) = (\ell(n-\ell)^{-1})N^{-1}\sum_{i=1}^{N}(\tilde{\gamma}_n^{(i)} - \hat{\gamma}_n)^2,$$

where $\tilde{\gamma}_n^{(i)} = \ell^{-1}(n\hat{\gamma}_n - (n-\ell)\hat{\gamma}_n^{(i)})$ is the *$i$th pseudo-value* of $\hat{\gamma}_n$. Note that, if $\ell \equiv 1$, that is, if we delete a single observation at a time, the variance estimator of (2.5) reduces to its i.i.d. counterpart,

$$(2.6) \qquad \widehat{\mathrm{Var}}_{\mathrm{J}}(\hat{\gamma}_n) = (n(n-1))^{-1}\sum_{i=1}^{n}(\tilde{\gamma}_n^{(i)} - \hat{\gamma}_n)^2.$$

Now we are ready to define the JAB estimator of the variance of a block bootstrap estimator. Here, the jackknife pseudo-values for the JAB method are obtained by deleting "blocks of bootstrap blocks" at a time instead of deleting "blocks of original data values." This simple modification is important, as it results in considerable computational efficiency (see [27] for details). Specifically, let $\hat{\varphi}_n$ be the block-bootstrap estimator of $\varphi_n \equiv \varphi(G_n)$ based on blocks of size $\ell$, where $G_n(\cdot) = P(T_n \leq \cdot)$ is the distribution function of $T_n$ and $\varphi$ is a functional, such as the quantile functional. Let $m$ be an integer between 1 and $N$, and let $M = N - m + 1$. We define the JAB pseudo-values by deleting $m$ consecutive blocks from the observed blocks $\{\mathcal{B}_1, \ldots, \mathcal{B}_N\}$. For $i = 1, \ldots, M$, define the set $I_i = \{1, \ldots, N\} \setminus \{i, \ldots, i +$ $m-1\}$. To define the $i$th jackknife block-deleted point value $\hat{\varphi}_n^{(i)}$, one needs to resample $b$ blocks from the reduced collection $\{\mathcal{B}_j : j \in I_i\}$ and apply the functional $\varphi$ to the resulting block-bootstrap estimator $\hat{G}_{n,i}$ of $G_n$. More precisely, let $J_{i1}, \ldots, J_{ib}$ be a collection of $b$ random variables such that, conditional on $\{X_n\}_{n\geq 1}$, they are i.i.d. with the discrete uniform distribution on $I_i$, that is,

$$P_*(J_{i1} = j) = (N - m)^{-1} \qquad \text{for all } j \in I_i.$$

Then the desired resampled blocks are given by $\{\mathcal{B}_j^{*i} \equiv \mathcal{B}_{J_{ij}} : j = 1, \ldots, b\}$. Let $T_n^{*i}$ denote the bootstrap version of $T_n$ based on the resampled values in the blocks $\{\mathcal{B}_j^{*i} : j = 1, \ldots, b\}$. Let $\hat{G}_{n,i}$ be the conditional distribution of $T_n^{*i}$ (given $\{X_n\}_{n\geq 1}$). Then, the $i$th jackknife block-deleted point value $\hat{\varphi}_n^{(i)}$ is given by

$$\hat{\varphi}_n^{(i)} = \varphi(\hat{G}_{n,i}), \qquad i = 1, \ldots, M.$$

The JAB estimator of the variance of $\hat{\varphi}_n$ is defined as [cf. (2.5)]

$$(2.7) \qquad \widehat{\text{Var}}_{\text{JAB}}(\hat{\varphi}_n) = (m(N-m)^{-1})M^{-1}\sum_{i=1}^{M}(\tilde{\varphi}_n^{(i)} - \hat{\varphi}_n)^2,$$

where $\tilde{\varphi}_n^{(i)} = m^{-1}(N\hat{\varphi}_n - (N-m)\hat{\varphi}_n^{(i)})$. Thus, in the dependent case the JAB variance estimators of the bootstrap distribution function $\hat{\varphi}_{1n}$ and of the bootstrap quantile $\hat{\varphi}_{2n}$ are defined using formula (2.7). In the independent case the corresponding JAB estimators are given by (2.7) by setting $m \equiv 1 \equiv \ell$ for all $n \geq 1$ (cf. [11]).

In the next section we show that, under some regularity conditions, the JAB variance estimators are consistent.

**3. Main results.** In a seminal paper Götze and Hipp [15] obtained Edgeworth expansions for sums of dependent random vectors under some fairly general conditions. For the dependent case, we shall establish the consistency of $\widehat{\text{Var}}_{\text{JAB}}(\hat{\varphi}_{kn})$, $k = 1, 2$, under a similar framework. To this end, suppose that $X_1, X_2, \ldots$ are defined on a common probability space $(\Omega, \mathcal{F}, P)$ and a sequence $\mathcal{D}_0, \mathcal{D}_{\pm 1}, \mathcal{D}_{\pm 2}, \ldots$ of sub-$\sigma$-fields of $\mathcal{F}$ is given. Let $|x| = |x_1| + \cdots + |x_d|$ and $\|x\| = (x_1^2 + \cdots + x_d^2)^{1/2}$, respectively, denote the $\ell^1$ and $\ell^2$ norms of $x \equiv (x_1, \ldots, x_d)' \in \mathbb{R}^d$. Also, for a vector $\nu = (\nu_1, \ldots, \nu_d)' \in \mathbf{Z}_+^d$, let $D^\nu$ denote the $\nu$th order partial derivative operator $\partial^{\nu_1 + \cdots + \nu_d}/\partial x_1^{\nu_1} \cdots \partial x_d^{\nu_d}$, where $\mathbf{Z}_+ = \{0, 1, 2, \ldots\}$.

3.1. *Assumptions.* We shall make use of the following assumptions to prove the results:

(A.1) $H:\mathbb{R}^d \to \mathbb{R}$ is four-times continuously differentiable and for all $\nu \in \mathbf{Z}_+^d$, $|\nu|=4$,
$$|D^\nu H(x)| \leq C(1+\|x\|^a), \qquad x \in \mathbb{R}^d,$$
for some integer $a \geq 1$ and $C \in (0,\infty)$.

(A.2) There exists $\delta > 0$ such that, for all $n,k = 1,2,\ldots$ with $k > \delta^{-1}$, there exists a $\mathcal{D}_{n-k}^{n+k}$-measurable $d$-variate random vector $\tilde{X}_{n,k}$ such that
$$E\|X_n - \tilde{X}_{n,k}\| \leq \delta^{-1}\exp(-\delta k).$$

(A.3) There exists $\delta > 0$ such that, for all $n,k = 1,2,\ldots, A \in \mathcal{D}_{-\infty}^n, B \in \mathcal{D}_{n+k}^\infty$,
$$|P(A \cap B) - P(A)P(B)| \leq \delta^{-1}\exp(-\delta k).$$

(A.4) There exists $\delta > 0$ such that, for all $n,k,r = 1,2,\ldots$ and $A \in \mathcal{D}_{n-r}^{n+r}$,
$$E|P(A|\mathcal{D}_j : j \neq n) - P(A|\mathcal{D}_j : 0 < |n-j| \leq k+r)| \leq \delta^{-1}\exp(-\delta k).$$

(A.5) There exists $\delta > 0$ such that, for all $n,k = 1,2,\ldots, \delta^{-1} < k < n$, and for all $t \in \mathbb{R}^p$ with $\|t\| > \delta$,
$$E|E(\exp(\sqrt{-1}t'(X_{n-k} + \cdots + X_{n+k}))|\mathcal{D}_j : j \neq n)| \leq e^{-\delta}.$$

(A.6)$_q$ $E\|X_1\|^{2q+\delta} < \infty$ for some $\delta > 0$, where $q \in \mathbf{Z}_+$.

Formulation of assumptions (A.2)–(A.5) in terms of the $\sigma$-fields $\mathcal{D}_j$, instead of the random vectors $X_j$ themselves, is to provide more flexibility in verification of the assumptions and to extend the range of applicability of the results. Choosing $\mathcal{D}_j = \sigma\langle X_j\rangle$, the $\sigma$-field generated by $X_j$, renders the verification of (A.2) rather trivial, but makes the verification of the conditional Cramér condition (A.5) quite difficult. In a specific problem, a different choice of $\mathcal{D}_j$ (from $\sigma\langle X_j\rangle$) is often more suitable for verifying (A.2)–(A.5). For example, if $X_j = \sum_{k=-\infty}^\infty c_k Y_{j-k}$ for some sequence of i.i.d. random vectors $\{Y_j\}$ and constants $\{c_j\}$, then a natural choice of $\mathcal{D}_j$ is $\sigma\langle Y_j\rangle$, rather than $\sigma\langle X_j\rangle$. We refer the interested reader to [4] and [15, 16], where conditions (A.2)–(A.5) are verified for a number of important dependent models.

Next we comment on the moment condition (A.6)$_q$, which is used with different values of $q$ in the results below. For the distribution function case, we use (A.6)$_q$ primarily for establishing valid expansions for the variance of the bootstrap distribution function estimator $\hat{\varphi}_{1n}$. Essentially this requires establishing uniform integrability of the *reciprocal* of $\tilde{\tau}_n^k$ and of $E_*(\|\bar{X}_n^*\|^{5p})$ for suitable integers $k \geq 1$, $p \geq 1$, where $\tilde{\tau}_n^2 = n_1 \text{Var}_*(h(\hat{\mu})'\bar{X}_n^*)$. Our conditions are comparable to the moment conditions used in similar contexts (cf. [2]). For the quantile case, we need a stronger moment condition than for the distribution function case. Note that $\hat{\varphi}_{1n}$, being the probability of an event, is always bounded by 1, while the bootstrap quantile $\hat{\varphi}_{2n}$ can be very large when either the numerator of $T_n^*$ is large or its denominator is small. The higher moment condition for the quantile case results from making the $L^2$-norm of $\hat{\varphi}_{2n}$ over certain low-probability pathological sets small. The moment condition can be somewhat reduced if we restrict the range of possible values of the block size $\ell$ and the jackknife-variable $m$ in Theorems 3.1 and 3.3. However, we do not pursue such refinements here.

3.2. *Distribution function.* The first result deals with the JAB variance estimator of the block bootstrap distribution function $\hat{G}_n(x)$ of $T_n^*$ at a point $x = x_0 \in \mathbb{R}$.

THEOREM 3.1. *Suppose that assumptions* (A.1)–(A.5) *and* (A.6)$_q$ *hold with* $q = 31$. *Also, suppose that* $n^\kappa \leq \ell \leq \kappa^{-1} n^{1/3}$ *and that* $m^{-1}\ell^{1+\kappa} + mn^{-2/3} = o(1)$ *as* $n \to \infty$ *for some* $\kappa > 0$. *Then, for any* $x_0 \in \mathbb{R}$,

$$(3.1) \qquad \widehat{\mathrm{Var}}_{\mathrm{JAB}}(\hat{\varphi}_{1n})/\mathrm{Var}(\hat{\varphi}_{1n}) \xrightarrow{p} 1 \qquad \text{as } n \to \infty,$$

*where* $\hat{\varphi}_{1n} = P_*(T_n^* \leq x_0)$ *and* $T_n^*$ *is as defined in* (2.2).

Thus, Theorem 3.1 shows that, for a wide range of values of the block size $\ell$ and "blocks of blocks" jackknife variable $m$, the JAB variance estimator is consistent for the MBB distribution function estimator $\hat{\varphi}_{1n} = P_*(T_n^* \leq x_0)$. In particular, the JAB variance estimator of $\hat{\varphi}_{1n}$ is consistent when the MBB block size $\ell$ goes to infinity at the optimal rate $Cn^{1/4}$, $C \in (0, \infty)$ (cf. [20]).

A consistency result similar to Theorem 3.1 also holds for Efron's [11] JAB variance estimator in the i.i.d. case. To state it, let $X_{jk}$ denote the $k$th component of $X_j$, $j \geq 1, 1 \leq k \leq d$. Let $V_j$ be a $d + [d(d+1)/2]$-dimensional random vector with the first $d$ components given by $X_j$ and the remaining $[d(d+1)/2]$ components given by $\{X_{jp}X_{jq} : 1 \leq p \leq q \leq d\}$. Then we have the following result.

THEOREM 3.2. *Assume that* $X_1, X_2, \ldots$ *are i.i.d. random vectors in* $\mathbb{R}^d$ *satisfying the following version of the Cramér condition:*

$$(3.2) \qquad \limsup_{\|t\| \to \infty} |E \exp(it'V_1)| < 1.$$

*If, in addition, assumptions* (A.1) *and* (A.6)$_q$ *hold with* $q = 20$, *then for any* $x_0 \in \mathbb{R}$,

$$\widehat{\mathrm{Var}}_{\mathrm{JAB}}(\hat{\varphi}_{1n})/\mathrm{Var}(\hat{\varphi}_{1n}) \xrightarrow{p} 1$$

*as* $n \to \infty$, *where* $\hat{\varphi}_{1n} = P_*(T_n^* \leq x_0)$, *and* $T_n^*$ *and* $\widehat{\mathrm{Var}}_{\mathrm{JAB}}(\cdot)$ *are, respectively, defined by* (2.2) *and* (2.7) *with* $\ell \equiv 1$ *and* $m \equiv 1$ *for all* $n \geq 1$.

Thus, the original version of the JAB variance estimator, proposed by Efron [11] for independent r.v.s, is also consistent for the distribution function estimator obtained by Efron's [9] bootstrap method. Note that in the independent case the variance of the bootstrap distribution function estimator $\hat{\varphi}_{1n}$ goes to zero at a faster rate, namely $O(n^{-2})$, than that of the block bootstrap estimator. Therefore, the JAB variance estimator adapts itself to this higher accuracy level of the bootstrap estimator in the independent case.

A comparison of the regularity conditions for Theorems 3.1 and 3.2 shows that the moment condition used for the independent case is somewhat weaker than that in the dependent case. To see why, observe that, for the bootstrap method of Efron [9], the second- and higher-order terms in the Edgeworth expansions of the bootstrap estimator $\hat{\varphi}_{1n}$ and its jackknife pseudo-values are smooth functions of certain sample means. However, in the dependent case these terms in the Edgeworth expansions are given by functions of "sums of block-sums" which, due to their more complicated forms, lead to the higher moment requirement. Next consider the Cramér conditions for Theorems 3.1 and 3.2. In the independent case we need to assume the Cramér condition (3.2) on a second degree polynomial $V_1$ of $X_1$. This is because the stochastic expansion for the Studentized statistic $T_n^*$ is given by a smooth function of a similar quadratic function $V_j^*$ of the resampled bootstrap variables $X_j^*$. However, in the dependent case the corresponding stochastic expansion is given by a function of certain quadratic forms in the resampled *block sums*. Since the block length tends to infinity, the block sums approximately behave like a normal random vector, and an analog of the Cramér condition (3.2) is essentially ensured by (A.5). As a result, for the dependent case we do not need to make such an assumption on the quadratic function of the $X_j$'s separately.

3.3. *Bootstrap quantiles.* In this section we consider the JAB estimator of the variance of a bootstrap quantile. Theorem 3.3 below asserts the consistency of the JAB variance estimator of the block bootstrap quantile $\hat{\varphi}_{2n}$ in the dependent case.

THEOREM 3.3. *Suppose that assumptions* (A.1)–(A.5) *and* (A.6)$_q$ *hold with* $q = \max\{60, 12(2+a)\} + 1$, *where $a$ is as in* (A.1). *Also, suppose that $n^\kappa \leq \ell \leq \kappa^{-1} n^{1/3}$ and that $m^{-1} \ell^{1+\kappa} + mn^{-2/3} = o(1)$ as $n \to \infty$ for some $\kappa > 0$. Then, for any $\alpha \in (0,1)$,*

$$(3.3) \qquad \widehat{\mathrm{Var}}_{\mathrm{JAB}}(\hat{\varphi}_{2n})/\mathrm{Var}(\hat{\varphi}_{2n}) \xrightarrow{p} 1 \qquad as\ n \to \infty,$$

*where $\hat{\varphi}_{2n} = \inf\{x \in \mathbb{R} : P_*(T_n^* \leq x) \geq \alpha\}$ and $T_n^*$ is as defined in* (2.2).

Thus, it follows that, under fairly general conditions, the JAB estimator of the variance of an MBB quantile is consistent. In the literature on jackknife for independent r.v.s, a modification of the standard jackknife method, called the "delete-d" jackknife, is often used to produce a consistent estimator of the variance of the sample quantile. However, the consistency of the JAB variance estimator of a block bootstrap quantile is not merely an artifact of deleting blocks of blocks, as follows by considering the i.i.d. case. Indeed, for independent data, the jackknife part of the JAB method deletes only one observation at a time and the consistency of the resulting variance estimator for a bootstrap estimator continues to hold, as shown by the following theorem.

THEOREM 3.4. *Assume that $X_1, X_2, \ldots$ are independent $\mathbb{R}^d$-valued r.v.s satisfying the conditions of Theorem 3.2 and that $E\|X_1\|^p < \infty$ for $p = \max\{75, 4(4+2a)\}$. Let $\hat{\varphi}_{2n}$ and $\mathrm{Var}_{\mathrm{JAB}}(\hat{\varphi}_{2n})$ be, respectively, defined by (2.4) and (2.7) with $\ell \equiv 1$ and $m \equiv 1$ for all $n \geq 1$. Then (3.3) holds.*

Thus, under the regularity conditions of Theorem 3.4, the JAB variance estimator of a bootstrap quantile is consistent. As it is well known, the standard (delete-1) jackknife method fails drastically in variance estimation problems involving the *sample* quantile and other nonsmooth functionals of the empirical distribution function. However, the JAB variance estimator, which also applies the delete-1 jackknife method to a *bootstrap* quantile estimator under independence, effectively overcomes this inconsistency problem. This appears to be a striking property of the JAB method and highlights the marked contrast between the properties of the JAB method and the standard jackknife method.

REMARK 3.1. A version of the JAB method for the nonoverlapping block bootstrap (NBB) method (cf. [7, 28]) was described in [27]. Consistency of the JAB variance estimators for the NBB distribution function estimators and the NBB quantile estimators continues to hold solely under the conditions of Theorems 3.1 and 3.3, respectively. A similar remark also applies to the JAB variance estimators of bootstrap estimators based on the Circular Block Bootstrap method of Politis and Romano [35].

REMARK 3.2. In many applications involving dependent data, block bootstrap estimators are used with data-based choices of $\ell$. Although the consistency results proved in here are not directly applicable to such cases, the following observation may be worth noting. Suppose that the JAB variance estimator is applied to a bootstrap estimator with an estimated block size $\hat{\ell}$ and that $\hat{\ell}$ takes values in an interval of the form

$$L \equiv [C_1 n^{1/4}, C_2 n^{1/4}],$$

with probability one, for some positive constants $C_1, C_2 \in (0,1)$. (Here $n^{1/4}$ is the optimal order of the block size for estimating $\varphi_{1n}, \varphi_{2n}$.) In this case one can refine the arguments of the paper to establish *uniform* versions of the consistency results of Theorems 3.1 and 3.3 over the interval $L$. As a result, when the JAB method is employed with such estimators of the block size, the corresponding JAB variance estimator is expected to estimate the additional variability in the bootstrapped estimators $\hat{\varphi}_{jn}(\hat{\ell})$, $j=1,2$, consistently.

**4. Numerical results.** In this section we report the results from a simulation study conducted to investigate finite sample properties of the JAB method and the choice of the JAB smoothing parameter $m$ in the dependent case. For results under the independent case, see Efron [11], Hill, Cartwright and Arbaugh [21] and the references therein. For the simulation study, we considered four different time series models, given by:

(I)  $X_t = (\varepsilon_t + \varepsilon_{t+1})/\sqrt{2}$,
(II) $X_t = 0.3 X_{t-1} + \varepsilon_t$,
(III) $X_t = -0.1 X_{t-1} + \varepsilon_t$,
(IV) $X_t = W_t^2 + W_t \mathbb{1}(W_t < 0)$, with $W_t = 0.6 W_{t-1} - 0.3 W_{t-2} + 0.1 W_{t-3} + \varepsilon_t + 0.2 \varepsilon_{t-1} + 0.3 \varepsilon_{t-2} + 0.1 \varepsilon_{t-3}$,

where, in (I)–(IV), $\{\varepsilon_t\}$ is a sequence of i.i.d. $N(0,1)$ random variables. Thus, model (I) is a moving average (MA) model of order 1, models II and III are autoregressive (AR) models of order 1, with a positive and a negative autoregression parameter, respectively, and model (IV) is an instantaneous transformation of an ARMA(3,3) model. The marginal distribution of the $X_t$'s under the first three models is *symmetric*, but under model (IV), it is *asymmetric* (with a thicker right tail). For each of the models (I)–(IV), we considered two sample sizes $n = 125, 800$ and applied the MBB with the block lengths $\ell = 5$ and $\ell = 8$, respectively, to estimate three population parameters: (i) $\varphi_{1n} = P(T_n \leq 0)$, the distribution function of $T_n$ at 0, (ii) $\varphi_{2n}$ = the 0.35 quantile of $T_n$ and (iii) $\varphi_{3n}$ = the 0.80 quantile of $T_n$. The corresponding MBB estimators are given by

$$\check{\varphi}_{1n} = P_*(T_n^* \leq 0),$$
(4.1) $$\check{\varphi}_{2n} = 0.35 \text{ quantile of } T_n^*,$$
$$\check{\varphi}_{3n} = 0.80 \text{ quantile of } T_n^*,$$

where $T_n^*$ is as in (2.2). Here the *parameters of interest* for the JAB method are given by $\text{Var}(\check{\varphi}_{kn})$, $k=1,2,3$. Table 1 gives the true values of $\text{Var}(\check{\varphi}_{kn})$ obtained by 20,000 simulation runs for each model. For each simulated data set, $K = 1000$ bootstrap replicates were used to compute the bootstrap estimators $\check{\varphi}_{kn}$. Next we applied the JAB method to estimate the variances of the MBB estimators $\check{\varphi}_{kn}$. As in [27], we considered several values of the jackknife block length $m$ using the formula

$$m = Cn^{1/3}\ell^{2/3} \tag{4.2}$$

for different constants $C$. An intuitive explanation for the choice of $m$ in (4.2) may be given as follows. Note that the JAB applies a version of the moving block jackknife method to the MBB estimator $\check{\varphi}_{kn}$ by deleting "blocks" of $m$ *overlapping* blocks. Since (i) the MSE-optimal block length for the block jackknife variance estimator (cf. [20, 22]) for a sample of size $J$ is of the form $C \cdot J^{1/3}$ and (ii) there are $b = \lfloor n/\ell \rfloor$ many "approximately independent" nonoverlapping blocks of length $\ell$ in a sample of size $n$ for a block bootstrap estimator based on *nonoverlapping* blocks of length $\ell$, the optimal jackknife block length for *its* jackknife variance estimator is of the form $C \cdot b^{1/3}$. But there are approximately $C \cdot b^{1/3}\ell$ *overlapping* blocks contained in $C \cdot b^{1/3}$ *nonoverlapping* blocks and vice versa. Hence, for a MBB estimator based on *overlapping* blocks of length $\ell$, one should delete $C \cdot b^{1/3}\ell$ overlapping blocks to define the JAB variance estimator. Formula (4.2) now follows by noting that $C \cdot b^{1/3}\ell$ is asymptotically equivalent to $Cn^{1/3}\ell^{2/3}$.

We applied the JAB method with the above choice of $m$ for $C = 0.1, 0.5, 1.0$. Table 2 gives the bias, the standard deviation (sd) and the mean squared error (mse) of the JAB estimators of $\text{Var}(\check{\varphi}_{kn})$ under models (I)–(IV). The results are based on 1000 simulation runs in each case. For each simulation run the number of bootstrap replicates was $K = 1000$. We make the following observations on the results of Table 2.

(1) For all four models and for all three parameters $\text{Var}(\check{\varphi}_{kn})$, $k = 1, 2, 3$, the mse's of the JAB estimators decreased as the value of $C$ is decreased. The best performance is attained at $C = 0.1$ in each case.

(2) There seems to be little effect of asymmetry of the marginal distribution of $X_t$ on the choice of $C$ in (4.2).

TABLE 1
*The true values of the target parameters* $\text{Var}(\check{\varphi}_{kn}), k = 1, 2, 3$, *for models* (I)–(IV) *based on* 20,000 *simulation runs. The MBB estimators* $\check{\varphi}_{kn}$'s *are defined by* (4.1) *with* $(n, \ell) = (125, 5)$ *and* $(n, \ell) = (800, 8)$. *For each data set,* $K = 1000$ *bootstrap replicates were used*

| | $n = 125, \ell = 5$ | | | $n = 800, \ell = 8$ | | |
|---|---|---|---|---|---|---|
| Model | $\text{Var}(\check{\varphi}_{1n})$ | $\text{Var}(\check{\varphi}_{2n})$ | $\text{Var}(\check{\varphi}_{3n})$ | $\text{Var}(\check{\varphi}_{1n})$ | $\text{Var}(\check{\varphi}_{2n})$ | $\text{Var}(\check{\varphi}_{3n})$ |
| (I) | 2.051e–4 | 2.219e–3 | 1.956e–3 | 1.715e–4 | 1.226e–3 | 1.363e–3 |
| (II) | 1.846e–4 | 1.351e–3 | 2.292e–3 | 1.688e–4 | 1.232e–3 | 1.345e–3 |
| (III) | 1.837e–4 | 1.297e–3 | 2.098e–3 | 1.708e–4 | 1.240e–3 | 1.350e–3 |
| (IV) | 1.870e–4 | 1.348e–3 | 2.272e–3 | 1.748e–4 | 1.295e–3 | 1.377e–3 |

TABLE 2
*Bias, standard deviation (sd) and mean squared error (mse) of the JAB variance estimators of the MBB estimators $\check{\varphi}_{kn}, k = 1, 2, 3$, of (4.1) under models (I)–(IV), where the jackknife parameter m is chosen according to the formula $m = Cn^{1/3}\ell^{2/3}$ with different values of C. The results are based on 1000 simulation runs in each case. The value of K is the same as in Table 1. An asterisk "∗" indicates the minimum value of the mse*

|  |  | $n = 125, \ell = 5$ |  |  | $n = 800, \ell = 8$ |  |  |
| --- | --- | --- | --- | --- | --- | --- | --- |
|  |  | $C = 0.1$ | $C = 0.5$ | $C = 1.0$ | $C = 0.1$ | $C = 0.5$ | $C = 1.0$ |
| Model (I) | | | | | | | |
| $\mathrm{Var}(\check{\varphi}_{1n})$ | bias | 1.235e–4 | 8.963e–3 | 2.265e–2 | 2.856e–4 | 3.626e–3 | 7.520e–3 |
|  | sd | 1.421e–5 | 1.427e–3 | 7.205e–3 | 1.993e–5 | 4.891e–4 | 1.606e–3 |
|  | mse | 1.546e–8∗ | 8.237e–5 | 5.648e–4 | 3.825e–6∗ | 6.276e–4 | 2.902e–3 |
| $\mathrm{Var}(\check{\varphi}_{2n})$ | bias | 2.407e–4 | 7.054e–2 | 0.1992 | 1.865e–3 | 2.476e–2 | 5.244e–2 |
|  | sd | 6.466e–4 | 1.988e–2 | 9.574e–2 | 5.873e–4 | 3.791e–3 | 1.233e–2 |
|  | mse | 4.760e–7∗ | 5.372e–3 | 4.885e–2 | 3.825e–6∗ | 6.276e–4 | 2.902e–3 |
| $\mathrm{Var}(\check{\varphi}_{3n})$ | bias | 5.825e–5 | 5.643e–2 | 0.1609 | 1.511e–3 | 2.215e–2 | 4.707e–2 |
|  | sd | 5.501e–4 | 1.173e–2 | 6.997e–2 | 6.017e–4 | 3.303e–3 | 1.085e–2 |
|  | mse | 3.060e–7∗ | 3.322e–3 | 3.078e–2 | 2.646e–6∗ | 5.016e–4 | 2.333e–3 |
| Model (II) | | | | | | | |
| $\mathrm{Var}(\check{\varphi}_{1n})$ | bias | 1.445e–4 | 8.984e–3 | 2.291e–2 | 2.897e–4 | 3.705e–3 | 7.770e–3 |
|  | sd | 1.564e–5 | 1.391e–3 | 7.425e–3 | 1.973e–5 | 5.016e–4 | 1.669e–3 |
|  | mse | 2.113e–8∗ | 8.264e–5 | 5.802e–4 | 8.435e–8∗ | 1.398e–5 | 6.316e–5 |
| $\mathrm{Var}(\check{\varphi}_{2n})$ | bias | 7.667e–4 | 6.030e–2 | 0.1706 | 1.848e–3 | 2.516e–2 | 5.435e–2 |
|  | sd | 5.093e–4 | 1.071e–2 | 6.569e–2 | 5.542e–4 | 3.885e–3 | 1.284e–2 |
|  | mse | 8.471e–7∗ | 3.751e–3 | 3.342e–2 | 3.724e–6∗ | 6.484e–4 | 3.119e–3 |
| $\mathrm{Var}(\check{\varphi}_{3n})$ | bias | 2.142e–5 | 6.168e–2 | 0.1750 | 1.583e–3 | 2.268e–2 | 4.869e–2 |
|  | sd | 6.714e–4 | 1.146e–2 | 6.560e–2 | 5.877e–4 | 3.424e–3 | 1.120e–2 |
|  | mse | 4.512e–7∗ | 3.936e–3 | 3.493e–2 | 2.851e–6∗ | 5.264e–4 | 2.496e–3 |
| Model (III) | | | | | | | |
| $\mathrm{Var}(\check{\varphi}_{1n})$ | bias | 1.447e–4 | 7.617e–3 | 1.883e–2 | 2.631e–4 | 3.393e–3 | 7.126e–3 |
|  | sd | 1.538e–5 | 1.294e–3 | 6.270e–3 | 2.098e–5 | 4.730e–4 | 1.542e–3 |
|  | mse | 2.117e–8∗ | 5.970e–5 | 3.939e–4 | 6.971e–8∗ | 1.174e–5 | 5.316e–5 |
| $\mathrm{Var}(\check{\varphi}_{2n})$ | bias | 8.643e–4 | 5.075e–2 | 0.1352 | 1.721e–3 | 2.298e–2 | 4.968e–2 |
|  | sd | 5.573e–4 | 9.881e–3 | 5.050e–2 | 5.738e–4 | 3.745e–3 | 1.214e–2 |
|  | mse | 1.058e–6∗ | 2.673e–3 | 2.083e–2 | 3.292e–6∗ | 5.421e–4 | 2.616e–3 |
| $\mathrm{Var}(\check{\varphi}_{3n})$ | bias | 1.988e–4 | 5.149e–2 | 0.1383 | 1.394e–3 | 2.058e–2 | 4.453e–2 |
|  | sd | 6.448e–4 | 1.032e–2 | 5.247e–2 | 5.874e–4 | 3.169e–3 | 1.040e–2 |
|  | mse | 4.553e–7∗ | 2.758e–3 | 2.188e–2 | 2.289e–6∗ | 4.336e–4 | 2.091e–3 |
| Model (IV) | | | | | | | |
| $\mathrm{Var}(\check{\varphi}_{1n})$ | bias | 1.415e–4 | 8.977e–3 | 2.341e–2 | 2.713e–4 | 3.528e–3 | 7.375e–3 |
|  | sd | 1.507e–5 | 1.382e–3 | 7.598e–3 | 1.929e–5 | 4.962e–4 | 1.595e–3 |
|  | mse | 2.025e–8∗ | 8.249e–5 | 6.056e–4 | 7.400e–8∗ | 1.269e–5 | 5.694e–5 |
| $\mathrm{Var}(\check{\varphi}_{2n})$ | bias | 8.130e–4 | 6.014e–2 | 0.1745 | 1.804e–3 | 2.451e–2 | 5.238e–2 |
|  | sd | 5.165e–4 | 1.054e–2 | 6.565e–2 | 6.083e–4 | 4.092e–3 | 1.272e–2 |
|  | mse | 9.277e–7∗ | 3.728e–3 | 3.476e–2 | 3.625e–6∗ | 6.177e–4 | 2.905e–3 |
| $\mathrm{Var}(\check{\varphi}_{3n})$ | bias | 7.782e–5 | 6.204e–2 | 0.1785 | 1.434e–3 | 2.125e–2 | 4.579e–2 |
|  | sd | 7.051e–4 | 1.177e–2 | 6.7749e–2 | 6.044e–4 | 3.377e–3 | 1.073e–2 |
|  | mse | 5.032e–7∗ | 3.988e–3 | 3.646e–2 | 2.424e–6∗ | 4.630e–4 | 2.212e–3 |

(3) As expected, there is indeed some effect of the arguments $x$ and $\alpha$ on the relative magnitudes of the mse's of the block bootstrap estimators $\hat{G}_n(x) = P_*(T_n^* \leq x)$ and $\hat{G}_n^{-1}(\alpha)$. For $C = 0.1$, $n = 125$ and $\ell = 5$, the mse of the JAB variance estimator of $\check{\varphi}_{1n}$ is the least, followed by those of $\check{\varphi}_{2n}$ and $\check{\varphi}_{3n}$, in all four models. Although it is uniform across the four models and possibly of secondary interest, a different pattern exists for $C = 0.5$ and $C = 1$.

(4) The choice $C = 0.1$ yields the best result for both combinations of the sample size and block length, $(n, \ell) = (125, 5)$ and $(800, 8)$. It is somewhat surprising that, despite having a larger sample size, the least mse values for the $(n, \ell) = (800, 8)$ case with $C = 0.1$ are consistently larger than those for the $(n, \ell) = (125, 5)$ case. This seems to be an artifact of the effects of the block size $\ell$ and the choice of $m$ through (4.2).

(5) When the results of this section are compared with the simulation results of Lahiri [27] on the choice of $m$ for the variance functional, a pronounced effect of the functional type on the value of $C$ is noticeable. It appears that, for the variance functional, a larger value of $C$ ($C \geq 1$) is desirable, whereas the value of $C$ should be small for the bootstrap distribution function and quantile estimators.

**5. Proofs.** We begin with some notation to be used in the rest of the paper. Let $\mathbb{1}_A$ and $\mathbb{1}(A)$ denote the indicator function of a set $A$, and for any countable set $J$ let $|J|$ denote the number of elements in $J$. Let $g : [0, \infty) \to [0, \infty)$ be an increasing function such that $g(t) = t$ if $t \leq 1, g(t) = 2$ for $t \geq 2$, and $g(\cdot)$ is Lipschitz on $[0, \infty)$. For $c > 0$, define the truncation function $g_c(\cdot) : \mathbb{R}^d \to \mathbb{R}^d$ by $g_c(x) = cxg(\|x\|/c)/\|x\|, x \in \mathbb{R}^d$. For $\alpha \in \mathbb{Z}_+^d$, set $\alpha! = \prod_{j=1}^d \alpha_j$ and $|\alpha| = \alpha_1 + \cdots + \alpha_d$. Let $\Phi$ and $\phi$, respectively, denote the distribution function and the density of a $N(0,1)$ random variable.

Next we collect the definitions of the random vectors that appear in the proofs below. For the sake of brevity, in this section we adopt the convention that, unless otherwise stated, the index $i$ runs over the set $\{0, 1, \ldots, M\}$, $j$ over $\{1, \ldots, b\}$ and $k$ over $\mathbb{N} \equiv \{1, 2, \ldots\}$. In the definitions below, we shall use the index $i = 0$ to define quantities that are related to the original block bootstrap variables (i.e., when no blocks are deleted). Thus, we set $I_0 = \{1, \ldots, N\}$ for the index set of all observed blocks and $\mathcal{B}_j^{*0} = \mathcal{B}_j^*$ for $j \in I_0$. With this convention, the resampled observations in the blocks $\{\mathcal{B}_j^{*i} : j = 1, \ldots, b\}$ are denoted by $X_1^{*i}, \ldots, X_{n_1}^{*i}$ (for all $i = 0, 1, \ldots, M$). Next define

$$U_k = (X_k + \cdots + X_{k+\ell-1})/\ell, \qquad U_{1k} = U_k \ell^{1/2},$$
$$U_j^{*i} = (X_j^{*i} + \cdots + X_{j+\ell-1}^{*i})/\ell, \qquad U_{1j}^{*i} = U_j^{*i} \ell^{1/2}.$$

Let $W_{2j}^{*i}$ (and $W_{2k}$) denote the $d(d+1)/2$-dimensional vector consisting of the diagonal and the subdiagonal elements of the $d \times d$ matrix $U_{1j}^{*i} U_{1j}^{*i'}$ (and of $U_{1k}U'_{1k}$, resp.). Next define $U_{2j}^{*i} = (U_{1j}^{*i\prime}; W_{2j}^{*i\prime})'$ and $U_{2k} = (U'_{1k}; W'_{2k})'$. Let $d_0 = d + d(d+1)/2$. For $t \in \mathbb{R}^{d_0}$, let $\hat{w}_{in}(t) = E_* \exp(\sqrt{-1}t'U_{21}^{*i})$ and $w_\ell(t) = E\exp(\sqrt{-1}t'U_{21})$. For $\alpha \in \mathbf{Z}_+^d$, set $c_\alpha = D^\alpha H(\mu)/\alpha!$ and $\hat{c}_{\alpha,i} = D^\alpha H(\hat{\mu}_{n,i})/\alpha!$. Let $\tilde{U}_j^{*i} = U_j^{*i} - \hat{\mu}_{n,i}$, $\tilde{U}_{rj}^{*i} = U_{rj}^{*i} - E_*U_{rj}^{*i}$, $r = 1, 2$.

Next define the covariance matrices

$$\hat{\Xi}_{n,i} = E_* \tilde{U}_{21}^{*i}\tilde{U}_{21}^{*i\prime}, \qquad \Xi_\ell = EU_{21}U'_{21}, \qquad \Xi_\infty = \lim_{\ell \to \infty} \Xi_\ell,$$

$$\tilde{\Sigma}_{n,i} = E_* \tilde{U}_{11}^{*i}\tilde{U}_{11}^{*i\prime} \quad \text{and} \quad \tilde{\tau}_{n,i}^2 = E_*(\tilde{Y}_{11}^{*i})^2,$$

where $\tilde{Y}_{1j}^{*i} = Y_{1j}^{*i} - E_*Y_{1j}^{*i}$ and $Y_{1j}^{*i} = \sum_{|\alpha|=1} \hat{c}_{\alpha,i}(U_{1j}^{*i})^\alpha$. Also, let $\tilde{Y}_{2j}^{*i} = (\tilde{Y}_{ij}^{*i})^2 - \tilde{\tau}_{n,i}^2$, $\bar{Y}_n^* = \sum_{|\alpha|=1} \hat{c}_{\alpha,i}(Z_n^{*i})^\alpha/\sqrt{n}_1$, $Z_n^{*i} = \sqrt{n}_1(\bar{X}_n^{*i} - \hat{\mu}_n)$ and $Z_{rn}^{*i} = b^{-1/2} \times \sum_{j=1}^b \tilde{U}_{rj}^{*i}$, $r = 1, 2$. Note that $Z_n^{*i} = Z_{1n}^{*i}$. Let $Z_n = n^{1/2}(\bar{X}_n - \mu)$ and let $Z_\infty$ denote a random vector with the $N(0, \Sigma_\infty)$ distribution on $\mathbb{R}^d$. Define the vector $Z_{2,\infty}$ from $Z_\infty$ in the same way $U_{21}$ is defined from $U_{11}$.

Next define the variables

$$\tilde{V}_n^{*i} \equiv \frac{\sqrt{n}_1}{b} \sum_{j=1}^b \{(\tilde{Y}_{ij}^{*i})^2 - \tilde{\tau}_{n,i}^2\} = \frac{\sqrt{n}_1}{b} \sum_{j=1}^b \tilde{Y}_{2j}^{*i}$$

and

$$\tilde{V}_{1n}^{*i} = \sum_{|\gamma|=1} [\hat{e}_{n,i}(\gamma)](Z_n^{*i})^\gamma,$$

where $\hat{e}_{n,i}(\gamma) = 2\sum_{|\alpha|=1} \sum_{|\beta|=1} \hat{c}_{\alpha,i}\hat{c}_{\beta+\gamma,i}(\beta+\gamma)!E_*(\tilde{U}_{11}^{*i})^\alpha(\tilde{U}_{11}^{*i})^\beta$. Let $\Gamma_{11}^{*i} = \sum_{|\alpha|=1} c_\alpha(U_{11}^{*i})^\alpha$, $\Gamma_{1k} = \sum_{|\alpha|=1} c_\alpha U_{1k}^\alpha$, $\bar{U}_{1n}^{(i)}(\alpha) = E_*(U_{11}^{*i})^\alpha$ and $\hat{U}_{1n}^{(i)} = E_*(U_{11}^{*i})^\alpha - EU_{11}^\alpha$. Also, write $\tilde{U}_{11}^{(i)}(\alpha) = [N\hat{U}_{11}^{(0)}(\alpha) - (N-m)\hat{U}_{11}^{(i)}]/m$ for the $i$th jackknife pseudo-value corresponding to $\hat{U}_{1n}^{(i)}$, $1 \le i \le M$. For *notational simplicity*, we often drop 0 from the subscripts and superscripts of the above variables for the case $i = 0$ and use both interchangeably, such as $U_j^* = U_j^{*0}$, $U_{1j}^* = U_{1j}^{*0}$ and $\hat{c}_{\alpha,0} = \hat{c}_\alpha$, and so on.

In the proofs below, unless otherwise stated, we set $\mu = 0$. We will use $C, C(\cdot)$ to denote generic constants, depending on their arguments (if any), but not on $i, n, \ell$ and $m$. For notational simplicity, any such dependence on population quantities, like $\|\Xi_\infty^{-1}\|$ or $E\|X_1\|^s$, that remain fixed throughout the proof of a result, will be suppressed. Similarly, let $\Lambda_j(x), j \ge 1$, denote generic polynomials (with numerical coefficients) in multivariables $x \in \mathbb{R}^{k_j}, 1 \le k_j < \infty$, that do not depend on $n, \ell, m$ and $i$. Unless otherwise specified, limits in order symbols are taken letting $n \to \infty$.

LEMMA 5.1. *Let $f: \mathbb{R}^d \to \mathbb{R}$ be a differentiable function such that* $\max\{|D^\nu f(x)| : |\nu| = 0, 1\} \le M(f)(1 + \|x\|^s)$, $x \in \mathbb{R}^d$, *for some $s \in [0, \infty)$*

and for some $M(f) \geq 1$. Assume that $4(\log n)^2 \leq \ell \leq n^{1/2}, \ell \leq m \leq n/4$, and that (A.2) and (A.3) hold. Then for any integer $r \geq 1$,

$$E(E_*f(U_{11}^{*i}) - Ef(U_{11}))^{2r} \leq C(r,s,d,\delta)[M(f)]^{2r}|I_i|^{-r}\ell^r$$

uniformly in $i \in \{0,1,\ldots,M\}$, provided $E\|X_1\|^{2q+\delta} < \infty$ for some $\delta > 0$ and for some integer $q$ satisfying $q > sr$ if $s > 0$ and $q \geq 0$ if $s = 0$.

PROOF. We prove the lemma for "$s > 0$" only. The case "$s = 0$" is simpler and can be proved using like arguments. For $c > 0$, let $f_c(\cdot) = f(g_c(\cdot))$. Then, with $c = n^r$, by Hölder's inequality,

$$E|(E_*f(U_{11}^{*i}) - Ef(U_{11}))^{2r} - (E_*f_c(U_{11}^{*i}) - Ef_c(U_{11}))^{2r}|$$
$$\leq C(r)E\{|(E_*f(U_{11}^{*i}) - E_*f_c(U_{11}^{*i})) - (Ef(U_{11}) - Ef_c(U_{11}))|$$
$$\times (|E_*f(U_{11}^{*i}) - Ef(U_{11})|^{2r-1} + |E_*f_c(U_{11}^{*i}) - Ef_c(U_{11})|^{2r-1})\}$$
$$\leq C(r)|I_i|^{-1}\sum_{j \in I_i}(M(f))^{2r}E\{(1+\|U_{11}\|^s)(1+\|U_{1j}\|^s)^{2r-1}\mathbb{1}(\|U_{11}\| > c)\}$$
$$+ C(r)(M(f))^{2r}\{E(1+\|U_{11}\|^s)\mathbb{1}(\|U_{11}\| > c)\} \cdot E(1+\|U_{11}\|^s)^{2r-1}$$
$$\leq C(r)(M(f))^{2r}c^{-\delta}E(1+\|U_{11}\|^s)^{2r+\delta}$$
$$\leq C(r,s,d)(M(f))^{2r}n^{-r}.$$

Next we approximate $X_j$'s with $\tilde{X}_{j,k}$'s. Let $\tilde{U}_{1j}$ and $\tilde{U}_{11}^{*i}$'s be defined by replacing $X_j$'s in the definitions of $U_{1j}$ and $U_{11}^{*i}$'s with $\tilde{X}_{j,u}, u = (\log n)^{3/2}$. Then, using (A.3), it is easy to check that

$$|E(E_*f_c(U_{11}^{*i}) - Ef_c(U_{11}))^{2r} - E(E_*f_c(\tilde{U}_{11}^{*i}) - Ef_c(\tilde{U}_{11}))^{2r}|$$
$$\leq C(r,s,d,\delta)[M(f)]^{2r}n^{-r}.$$

Let $K_i$ be the integer satisfying $(K_i - 1)2\ell < |I_i| \leq 2\ell K_i$. With $1 \leq k \leq K_i$, let $S(k,i)$ denote the sum of all $f_c(\tilde{U}_{1j})$ for $j \in [(k-1)2\ell+1, k2\ell] \cap I_i$. Noting that the $\tilde{X}_{j,u}$'s defining $S(k,j)$ and $S(k+2,j)$ are separated by a block of length $\ell$ or more and $u/\ell = o(1)$, for any $\eta > 0$ we get

$$E(E_*f_c(U_{11}^{*i}) - Ef_c(U_{11}))^{2r}$$
$$\leq C(r)|I_i|^{-2r}\left[E\left(\sum_1 [S(k,i) - ES(k,i)]\right)^{2r}\right.$$
$$\left.+ E\left(\sum_2 [S(k,i) - ES(k,i)]\right)^{2r}\right]$$
$$\leq C(r)|I_i|^{-2r}\max\{(E|S(k,i) - ES(k,i)|^{2r+\eta})^{2r/(2r+\eta)} : 1 \leq k \leq K_i\}$$

$$\times \left(1 + \sum_{j=1}^{\infty} j^{2r-1} \exp(-C(r,\eta)(\ell + (j-1)2\ell - 2u))\right)(|I_j|/\ell)^r$$

$$\leq C(r,\eta)|I_i|^{-r}\ell^r (E|f_c(\tilde{U}_{11})|^{2r+\eta})^{2r/(2r+\eta)}$$

$$\leq C(r,\eta)|I_i|^{-r}\ell^r [M(f)]^{2r}(1 + E\|g_c(\tilde{U}_{11})\|^{s(2r+\eta)})^{2r/(2r+\eta)},$$

where $\sum_1$ and $\sum_2$, respectively, denote summation over all odd and even $k$, $1 \leq k \leq K_i$. Now take $\eta = 2(q-rs)/s$, so that $s(2r+\eta) = 2q$. Then, by (A.3) and Lemma 3.30 of [15], the lemma follows. $\square$

LEMMA 5.2. *Suppose that assumptions* (A.2)–(A.4) *hold,* $E\|X_1\|^{4+\delta} < \infty$ *for some* $\delta > 0$ *and* $1 \leq \ell \leq m \leq n/4$.

(a) *For any integer* $r \geq 1$, *if* $E\|X_1\|^{2r+2+\delta} < \infty$ *for some* $\delta > 0$, *then*

$$\max_{0 \leq i \leq M} P(\|\hat{\mu}_{n,i} - \mu\| > n^{-1/2} \log n) = o(n^{-r} m^r).$$

(b) *For any* $\lambda > 0, K > 0$ *and* $\varepsilon > 0$, *there exist* $N_1 = N_1(\lambda, \varepsilon, K) \geq 1$ *and* $C(\lambda, \varepsilon) > 0$ *such that, for all* $n \geq N_1$ *and for any* $\ell \in [4(\log n)^2, n^{1/2}]$,

$$\max_{0 \leq i \leq M} P(\sup\{|\hat{w}_{in}(t) - w_\ell(t)| : \|t\| \leq \lambda n^\lambda\} > \varepsilon, E_*\|U_{21}^{*i}\| < K) \leq C(\lambda, \varepsilon)n^{-2\lambda}.$$

PROOF. (a) It is easy to check that, for any $0 \leq i \leq M$, $\|\hat{\mu}_{n,i} - \bar{X}_n\| \leq (N-m)^{-1}\{(\ell+m)\|\bar{X}_n\| + \mathbb{1}_{\{1,\ldots,M\}}(i)\|\sum_{j=i}^{i+m-1} U_j\| + \|\ell^{-1}\sum_{j=1}^{\ell-1}(\ell-j)(X_j + X_{n-j+1})\|\}$. By Theorem 2.4 of [24] (with $s = 2r + 2$),

$$P(\|\bar{X}_n\| > C(r)n^{-1/2}(\log n)^{1/2}) = o(n^{-r}).$$

Further, using arguments similar to the proof of Lemma 3.28 of [15], we get

$$E\left\|\sum_{j=1}^{m} U_j\right\|^{2r} \leq C(r)(\ell+m)^r$$

and

$$E\left\|\sum_{j=1}^{\ell-1}(\ell-j)(X_j + X_{n-j+1})\right\|^{2r} \leq C(r)\ell^{3r}.$$

Part (a) of the lemma now follows using Markov's inequality and the above bounds.

(b) Let $D_n = \{\nu \in \mathbf{Z}^d : \|\nu\| \leq \lambda n^{2\lambda}\}$. Then, for any $t \in \mathbb{R}^d$ with $\|t\| \leq \lambda n^\lambda$, there exists a $\nu \in D_n$ such that $\|t - \nu n^{-\lambda}\| < n^{-\lambda}$. Then it follows that for all $i \in \{0, 1, \ldots, M\}$,

$$|\hat{w}_{in}(t) - \hat{w}_{in}(\nu n^{-\lambda})| \leq n^{-\lambda} E_*\|U_{21}^{*i}\| \leq n^{-\lambda} K,$$
$$|w_\ell(t) - w_\ell(\nu n^{-\lambda})| \leq n^{-\lambda} E\|U_{21}\|.$$

Hence, by Lemma 5.1, for all $n \geq N(\lambda, \varepsilon, K) \geq 1$, and for any $\ell \in [4(\log n)^2, n^{1/2}]$, $m \in [\ell, n/4]$,

$$P(\sup\{|\hat{w}_{in}(t) - w_\ell(t)| : \|t\| \leq \lambda n^\lambda\} > \varepsilon, E_* \|U_{21}^{*i}\| < K)$$
$$\leq P(\max\{|\hat{w}_{in}(\nu n^{-\lambda}) - w_\ell(\nu n^{-\lambda})| : \nu \in D_n\} > \varepsilon/2)$$
$$\leq C(\varepsilon, r)|D_n| |I_i|^{-2r} \ell^{2r} \leq C(\varepsilon, \lambda, r) n^{2\lambda d}(n^{-1}\ell)^{2r}$$

uniformly in $i \in \{0, 1, \ldots, M\}$, where $r$ is an integer $\geq 2\lambda(d+1)$. This completes the proof. $\square$

LEMMA 5.3. *Let $W_1, W_2, \ldots, W_n$ be i.i.d. random vectors on $\mathbb{R}^d$ with $EW_1 = 0$, $\mathrm{Cov}(W_1) = I_d$ and $\rho_W = E\|W_1\|^5 < \infty$. Suppose that $n \geq 1$ is such that $a_{1n} \equiv E\|W_1\|^4 \mathbb{1}(\|W_1\| > 2\sqrt{n}/3) < (2^8 d)^{-1} n$. Then for any Borel set $B$ in $\mathbb{R}^d$ and any $\varepsilon > 0$,*

$$\left| P\left(n^{-1/2} \sum_{j=1}^n W_j \in B\right) - \int_B \left(1 + \sum_{k=1}^2 n^{-j/2} p_k(x)\right) d\Phi(x) \right|$$
(5.1)
$$\leq C(d)(1 + \rho_W)$$
$$\times [(n^{-3/2}) + (a_{1n}n^{-1} + a_{2n}^{2(d+5)} + n^{(2+d)}a_{2n}^{+(d+5)(d+8)})$$
$$+ (\alpha_n \varepsilon^{-2d} + n^{d+5}\varepsilon^{-8d} \exp(-\varepsilon^{-1})) + \Phi((\partial B)^{2\varepsilon})],$$

*where $a_{2n} = ((1 + \rho_W)n^{-1})^{1/(d+5)}$ and $\alpha_n = n^{d+5}\xi_n^{(n-d-5)}$ with $\xi_n = \sup\{|E \exp(it'W_1)| : (16E\|W_1\|^3)^{-1} \leq \|t\| \leq \varepsilon^{-4}\} + 2P(\|W_1\| > \sqrt{n})$. Here $p_1(x)\phi(x) = -\sum_{|\nu|=3}(EW_1^\nu/\nu!)(D^\nu \phi(x))$ and $p_2(x)\phi(x) = \sum_{|\nu|=4}\chi_\nu(W_1) \times D^\nu \phi(x)/\nu! + 2^{-1}(\sum_{|\nu|=3}(\nu!)^{-1}\chi_\nu(W_1)D^\nu)^2 \phi(x)$.*

PROOF. The lemma can be proved by carefully recasting the arguments and the bounds in the proof of Theorem 20.1 in [3] (cf. [2]). We omit the details. $\square$

LEMMA 5.4. *Assume that conditions (A.1)–(A.6) hold, $\rho(r,\delta) \equiv E\|X_1\|^{15(2r+1)+\delta} < \infty$ and $n^\kappa < \ell < \kappa^{-1} n^{1/3}$ and $\ell \leq m \leq n^{2/3}$ for some integer $r \geq 1$ and constants $0 < \kappa, \delta < \infty$. Then there exist a constant $C = C(\kappa, \delta, r) \in (0, \infty)$ and sets $A_{in}, i = 0, 1, \ldots, M$, such that, for all $n \geq C$,*

$$\max\{P(A_{in}^c) : i = 0, 1, \ldots, M\} \leq C[n^{-r}\ell^r],$$

*and on the set $A_{in}$,*

$$\max_{0 \leq i \leq M} \sup_x |P_*(T_n^{*i} \leq x) - \hat{\Psi}_{n,i}(x)| \leq C[b^{-3/2} + \ell n^{-3/2}](\log n)^4,$$

*where the Edgeworth expansion $\hat{\Psi}_{n,i}$ of $T_n^{*i}$ is defined by relation (5.14) below.*

PROOF. Without loss of generality, assume that the $X_j$'s and the bootstrap variables $X_1^{*i}, \ldots, X_{n_1}^{i*}$ are defined on a common probability space $(\Omega, \mathcal{F}, P)$ for all $n \geq 1$. For $i = 0, 1, \ldots, M$ define the set $A_{in}^{(1)}$ as the intersection of the sets $\{\|\tilde{\Sigma}_{n,i} - \Sigma_\ell\| \leq (2\|\Sigma_\ell^{-1}\|)^{-1}\}$, $\{|\tilde{\tau}_{n,i}^2 - \tau_\ell^2| \leq \tau_\ell^2/2\}$, $\{\|\hat{\mu}_{n,i} - \mu\| \leq \delta_1 n^{-1/2} \log n\}$, $\{E_*\|U_{11}^{*i}\|^{10} \leq \delta_1^{-1}\}$, $\{\sqrt{\ell}E_*(\tilde{Y}_{11}^*)^3| \leq \delta_1^{-1}\}$, and $\{\|\hat{\tilde{\Xi}}_{n,i} - \Xi_\ell\| < (2\|(\Xi_\ell)^{-1}\|)^{-1}\}$, where $0 < \delta_1 < 1$ is such that $E\|Z_\infty\|^{10} < (2\delta_1)^{-1}$, $\lim_{n \to \infty} \sum_{|\alpha|=3} |n^2 E \bar{X}_n^\alpha| < (2\delta_1)^{-1}$ and $\inf\{\sum_{|\alpha|=1} |D^\alpha H(x)| : \|x - \mu\| < 2\delta_1\} > \delta_1$. In the steps below, unless specifically mentioned otherwise, the bounds in the inequalities hold uniformly in $i \in \{0, 1, \ldots, M\}$.

STEP I (Stochastic expansion). First we derive a stochastic approximation to $T_n^{*i}$. For this, define the set $A_{in}^*$ as the intersection of the sets (involving both the $X_j$'s and the $X_j^{*i}$'s) $\{\|\bar{X}_n^{*i} - \hat{\mu}_{n,i}\| \leq \delta_1 n^{-1/2} \log n\}$, $\{|(\tau_n^{*i})^2 - \tilde{\tau}_\ell^2| \leq \tilde{\tau}_\ell^2/4\}$ and $\{|\tilde{V}_n^{*i}| \leq C(\delta_1) n^{1/2} b^{-1/2} \log n\}$. Then, using Taylor's expansion, we have

$$(5.2) \quad \sqrt{n}(H(\bar{X}_n^{*i}) - H(\hat{\mu}_{n,i})) = \sqrt{n} \sum_{|\alpha|=1}^{3} \hat{c}_{\alpha,i}(\bar{X}_n^{*i} - \hat{\mu}_{n,i})^\alpha + R_{1n}^{*i},$$

$$(\tau_n^{*i})^2 - \tilde{\tau}_{n,i}^2 = b^{-1} \sum_{j=1}^{b} [(\tilde{Y}_{1j}^{*i})^2 - \tilde{\tau}_{n,i}^2] + n_1^{-1/2} \tilde{V}_{1n}^{*i} - \ell(\bar{Y}_n^{*i})^2$$

$$(5.3) \qquad + n^{-1/2} b^{-1/2} \Lambda_1(Z_{2n}^{*i}; \{\hat{c}_{\alpha,i} : |\alpha| = 1\})$$

$$\qquad + n^{-1} \Lambda_2(Z_{2n}^{*i}; \{\hat{c}_{\alpha,i} : |\alpha| \leq 3\}) + R_{2n}^{*i},$$

where, on the set $A_{in}^* \cap A_{in}^{(1)}$, the remainder terms admit the bounds (uniformly in $i$)

$$(5.4) \qquad |R_{1n}^{*i}| \leq C(\delta_1, D_H) n^{-3/2} (\log n)^4,$$

$$(5.5) \qquad |R_{2n}^{*i}| \leq C(\delta_1, D_H) \ell n^{-3/2} (\log n)^3.$$

In (5.4) and (5.5) $D_H = \sup\{|D^\alpha H(x)| : |\alpha| \leq 3, \|x - \mu\| \leq 2\delta_1\} + \sup\{(\Sigma_{|\alpha|=1} |D^\alpha H(x)|)^{-1} : \|x - \mu\| \leq 2\delta_1\}$. Next, using Taylor's expansion of the function $g(x) = x^{-1/2}$ around $x = \tilde{\tau}_{n,i}^2$, and using (5.2) and (5.3), we get

$$\sqrt{n_1}(H(\bar{X}_n^{*i}) - H(\hat{\mu}_{n,i}))/(\tau_n^{*i} + n^{-1})$$

$$= \left[\sqrt{n_1}\bar{Y}_n^{*i} + n_1^{1/2} \sum_{|\alpha|=2}^{3} \hat{c}_{\alpha,i}(n_1^{-1/2} Z_n^{*i})^\alpha + R_{1n}^{*i}\right]$$

$$\times [\tilde{\tau}_{n,i}^{-1} - (2\tilde{\tau}_{n,i}^3)^{-1}\{(\tau_n^{*i})^2 - \tilde{\tau}_{n,i}^2\}$$

$$+ 3(8\tilde{\tau}_{n,i}^5)^{-1}\{(\tau_n^{*i})^2 - \tilde{\tau}_{n,i}^2\}^2 + R_{3n}^{*i}]$$

$$
\begin{aligned}
(5.6) \quad &= \sqrt{n}\bar{Y}_n^{*i}[\tilde{\tau}_{n,i}^{-1} - (2\tilde{\tau}_{n,i}^3)^{-1}\{n_1^{-1/2}(\tilde{V}_n^{*i} + \tilde{V}_{1n}^{*i}) - b^{-1}(\sqrt{n}\bar{Y}_n^{*i})^2\} \\
&\qquad\qquad\qquad + 3(8\tilde{\tau}_{n,i}^5)^{-1}n_1^{-1}\{\tilde{V}_n^{*i}\}^2] \\
&\quad + n_1^{-1/2}\sum_{|\alpha|=2}\hat{c}_{\alpha,i}(Z_n^{*i})^\alpha/\tilde{\tau}_{n,i} + n^{-1/2}b^{-1/2}\Lambda_3(Z_{2n}^{*i};B_{n,i}) \\
&\quad + n^{-1}\Lambda_4(Z_{2n}^{*i};B_{n,i}) + R_{4n}^{*i} \\
&\equiv T_{1n}^{*i} + R_{4n}^{*i}, \qquad \text{say,}
\end{aligned}
$$

where $B_{n,i} = \{\hat{c}_{\alpha,i} : 1 \leq |\alpha| \leq 3\} \cup \{\hat{\mu}_{n,i}\} \cup \{\tilde{\tau}_{n,i}^{-1}\}$, and by (4.3), (4.4) and (4.5), on the set $A_{in}^* \cap A_{in}^{(1)}$,

$$
\begin{aligned}
(5.7) \quad &|R_{3n}^{*i}| \leq C(\delta_1, D_H)|(\tau_n^{*i})^2 - \tilde{\tau}_{n,i}^2|^3 \leq C(\delta_1, D_H)b^{-3/2}(\log n)^3, \\
&|R_{4n}^{*i}| \leq C(\delta_1, D_H)[\ell n^{-3/2}(\log n)^4 + b^{-3/2}(\log n)^4].
\end{aligned}
$$

Here, $T_{1n}^{*i}$ gives the desired stochastic expansion for the Studentized statistic $T_n^{*i}$. Note that $T_{1n}^{*i}$ is a polynomial in the $d_0$-dimensional i.i.d. random vectors $W_j^{*i} = [\hat{\Xi}_{n,i}]^{-1/2}(U_{2j}^{*i} - E_*U_{2j}^{*i}), j = 1, \ldots, b$.

STEP II (Edgeworth expansion for $T_{1n}^{*i}$). Note that conditional on $\mathcal{X}_n$, $U_{21}^{*i}, \ldots, U_{2b}^{*i}$ are i.i.d. random vectors with covariance matrix $\hat{\Xi}_{n,i}$. Also, on the set $A_{in}^{(1)}$, $\hat{\Xi}_{n,i}^{-1}$ exists and $\|\hat{\Xi}_{n,i}^{-1}\| < 2\|\Xi_\ell^{-1}\|$. Hence, taking $W_j = W_j^{*i}$, $j = 1, \ldots, b$, and $\varepsilon = b^{-2}$ in Lemma 5.3, for all $n \geq N_1 \equiv N_1(\delta_1, \kappa)$ (not depending on $\ell$ or $i$) we find that, on the set $A_{in}^{(1)}$,

$$
\left| P_*\left(b^{-1/2}\sum_{j=1}^b W_j^{*i} \in B\right) - \int_B \left(1 + \sum_{k=1}^2 b^{-k/2}\hat{p}_{k,i}(x)\right)d\Phi(x)\right|
$$
$$
\leq C(\delta_1)[b^{-3/2} + \hat{\alpha}_{n,i}b^{4d} + \Phi((\partial B)^{2\varepsilon})]
$$

for all Borel sets $B$ in $\mathbb{R}^{d_0}$. Here $\hat{\alpha}_{n,i}$ and $\hat{p}_{k,i}$ are as in Lemma 5.3 under the given choice of $W_j$'s. Using the bounds on $\|\hat{\Xi}_{n,i}^{-1}\|$ and $E_*\|U_{11}^{*i}\|^{10}$ on $A_{in}^{(1)}$, for all $n \geq N_1$, on the set $A_{in}^{(1)}$, we get

$$
\begin{aligned}
\hat{\alpha}_{n,i} &\leq b^{d+5}[\sup\{|\hat{\omega}_{n,i}(t)| : C_1(\delta_1, \kappa) \leq \|t\| \leq C_2(\delta_1, \kappa)b^8\} \\
&\qquad\qquad + C(\delta_1, \kappa)b^{-5/2}]^{(b-d-5)} \\
&\equiv b^{d+5}[\hat{\xi}_{n,i} + C(\delta_1, \kappa)b^{-5/2}]^{(b-d-5)}, \qquad \text{say.}
\end{aligned}
$$

Let $\xi_\infty \equiv \sup\{|E\exp(\sqrt{-1}t'Z_{2,\infty})| : \|t\| \geq C_1(\delta_1, \kappa)\}$. Note that $0 < \xi_\infty < 1$. Next define the set $A_{in}^{(2)} \equiv \{\hat{\xi}_{n,i} < \xi_\infty + [1-\xi_\infty]/2\}$ and let $A_{in} = A_{in}^{(1)} \cap A_{in}^{(2)}$.

Then for all $n \geq N_2(\delta_1, \kappa, \xi_\infty)$, on the set $A_{in}$, we get

$$
(5.8) \quad \left| P_*\left( b^{-1/2} \sum_{j=1}^{b} W_j^{*i} \in B \right) - \int_B \left( 1 + \sum_{k=1}^{2} b^{-k/2} \hat{p}_{k,i}(x) \right) d\Phi(x) \right|
$$
$$
\leq C(\delta_1, \kappa, \xi_\infty)[b^{-3/2} + \Phi((\partial B)^{2\varepsilon})]
$$

for all Borel sets $B$ in $\mathbb{R}^{d_0}$. Next note that there exists a constant $C = C(\delta_1, \kappa) > 0$ such that, for all $n \geq C$, $\|[\hat{\Xi}_{n,i}]^{-1}\| + E_*\|U_{21}^{*i}\|^4 \leq C$ on the set $A_{in}$. Hence, by (5.8) we get the following moderate deviation inequalities: For all $n \geq 1$, uniformly in $i = 0, 1, \ldots, M$, on the set $A_{in}$,

$$
(5.9) \quad P_*\left( \left\| \sum_{j=1}^{b} \tilde{U}_{2j}^{*i} \right\| > Cb^{1/2} \log n \right) \leq Cb^{-3/2},
$$

$$
(5.10) \quad P_*(\|\bar{X}_n^{*i} - \hat{\mu}_{n,i}\| > Cn^{-1/2} \log n) \leq Cb^{-3/2},
$$

where the second inequality follows from the first by noting that $\sqrt{n_1}(\bar{X}_n^{*i} - \hat{\mu}_{n,i}) = b^{-1/2} \sum_{j=1}^{b} \tilde{U}_{1j}^{*i}$. Now, by (5.9) and (5.10) it follows that, for all $n \geq 1$, on the set $A_{in}$,

$$
(5.11) \quad P_*(A_{in}^{*c}) \leq Cb^{-3/2}.
$$

Let $\bar{W}_b^{*i} = b^{-1} \sum_{j=1}^{b} W_j^{*i}$. Note that the stochastic expansion $T_{1n}^{*i}$ in (5.6) can be written as $T_{1n}^{*i} = \sum_{k=0}^{2} b^{-k/2} q_k(\sqrt{b}\bar{W}_b^{*i}; \{\hat{\mu}_{n,i}, \hat{\Xi}_{n,i}\}) + n_1^{-1/2} q_3(\sqrt{b}\bar{W}_b^{*i}; \{\hat{\mu}_{n,i}, \hat{\Xi}_{n,i}\})$ for some suitable functions $q_k(\omega; \{x, \Xi\}), k = 0, 1, 2, 3$, where each $q_k(\omega; \{x, \Xi\})$ is a polynomial in $\omega \in \mathbb{R}^{d_0}$ with coefficients that are continuous functions of $(x, \Xi)$ in a neighborhood of $(\mu, \Xi_\infty)$. In particular, the functions $q_k$'s do *not* depend on $i$ and $n$. Since on $A_{in}$ $\sum_{|\alpha|=1} |D^\alpha H(\hat{\mu}_n + x)| > C(\delta_1)$ for $\|x\| < \delta_1$, by applying the transformation technique of Bhattacharya and Ghosh [1], and recasting their arguments on pages 434–435 for the last term above, one can show that, for all $n \geq N_3 \equiv N_3(\delta_1, \kappa, \xi_\infty)$ (not depending on $\ell$ and $i$), on the set $A_{in}$,

$$
(5.12) \quad \sup_{x \in \mathbb{R}} |P_*(T_{1n}^{*i} \leq x) - \hat{\Psi}_{n,i}(x)| \leq C(\delta_1, \kappa, D_H)[\ell n^{-3/2} + b^{-3/2}](\log n)^4,
$$

where the expansion $\hat{\Psi}_{n,i}$ is defined through the transformation of the Edgeworth expansion of $b^{-1/2} \sum_{j=1}^{b} W_j^{*i}$ in (5.8). This completes Step II.

STEP III (Form of $\hat{\Psi}_{n,i}$). We now identify the terms of $\hat{\Psi}_{n,i}$ by computing the approximate cumulants of $T_{1n}^{*i}$. The validity of the resulting expression follows by recasting the arguments of Bhattacharya and Ghosh [1] to the present setup. Let $\chi_*^{[r]}(S_n^*)$ denote the conditional cumulant of order $r \geq 1$ of a generic bootstrap variable $S_n^*$, defined as

$$(\sqrt{-1})^r \chi_*^{[r]}(S_n^*) = \frac{\partial^r}{\partial t^r} \log E_* \exp(\sqrt{-1}t' S_n^*)\Big|_{t=0}.$$

Next note that, for any zero mean i.i.d. random vectors $W_j = (W_{1j}, \ldots, W_{6j})'$, $j = 1, \ldots, J$,

$$E \prod_{k=1}^{3}\left(\sum_{j=1}^{J} W_{kj}\right) = JE\left(\prod_{k=1}^{3} W_{k1}\right),$$

$$E \prod_{k=1}^{4}\left(\sum_{j=1}^{J} W_{kj}\right) = JE\left(\prod_{k=1}^{4} W_{kj}\right) + J(J-1) \sum_{(4)}^{*} EW_1(I)EW_1(I^c),$$

(5.13) $\quad E \prod_{k=1}^{5}\left(\sum_{j=1}^{J} W_{kj}\right) = JE\left(\prod_{k=1}^{5} W_{k1}\right) + J(J-1) \sum_{(5)}^{*} EW_1(I)EW_1(I^c),$

$$E \prod_{k=1}^{6}\left(\sum_{j=1}^{J} W_{kj}\right) = JE\left(\prod_{k=1}^{5} W_{k1}\right) + J(J-1) \sum_{(6)}^{*} EW_1(I)EW_1(I^c)$$

$$+ J(J-1)(J-2) \sum_{(6)}^{**} EW_1(I_1)EW_1(I_2)EW_1(I_3),$$

where for $a \geq 4$ $\sum_{(a)}^{*}$ extends over all *distinct* partitions $\{I, I^c\}$ of $\{1, \ldots, a\}$ with $|I| = 2$, $\sum_{(a)}^{**}$ extends over all *distinct* partitions $\{I_1, I_2, I_3\}$ of $\{1, \ldots, a\}$ with $|I_1| = |I_2| = |I_3|$, and where for any set $I \subset \{1, \ldots, 6\}$ $W_1(I) = \prod_{k \in I} W_{k1}$. Now, using the identities in (5.13) and the bounds on the bootstrap moments on the set $A_{in}$, after some tedious algebra one can show that

$$\chi_*^{[1]}(T_{1n}^{*i}) = -\frac{1}{2\tilde{\tau}_{n,i}^3 \sqrt{n}}\left[\sqrt{\ell} E_*(\tilde{Y}_{11}^{*i})^3 + \sum_{|\gamma|=1} \hat{e}_{n,i}(\gamma) E_*\{\tilde{Y}_{11}^{*i}(\tilde{U}_{11}^{*i})^\gamma\}\right]$$

$$+ n^{-1/2} \sum_{|\alpha|=2} \hat{c}_{\alpha,i} E_*(\tilde{U}_{11}^{*i})^\alpha / \tilde{\tau}_{n,i} + \hat{Q}_{n,i}^{[1]}$$

$$\equiv n^{-1/2} \hat{\chi}_{n,i}^{[1]} + \hat{Q}_{n,i}^{[1]}, \qquad \text{say,}$$

$$\chi_*^{[2]}(T_{1n}^{*i}) = 1 + 3b^{-1} + \hat{Q}_{n,i}^{[2]},$$

$$\chi_*^{[3]}(T_{1n}^{*i}) = -\frac{1}{2\tilde{\tau}_{n,i}^3 \sqrt{n}}\left[7\sqrt{\ell} E_*(\tilde{Y}_{11}^{*i})^3 + 3 \sum_{|\gamma|=1} \hat{e}_n(\gamma) E_*\{\tilde{Y}_{11}^{*i}(\tilde{U}_{11}^{*i})^\gamma\}\right]$$

$$- \frac{3}{\sqrt{n}} \hat{\chi}_{n,i}^{[1]} + \hat{Q}_{n,i}^{[3]}$$

$$+ \frac{1}{\tilde{\tau}_{n,i}^3 \sqrt{n}} \sum_{|\alpha|=1} \sum_{|\beta|=1} \hat{c}_{\alpha+\beta,i} [\tilde{\tau}_{n,i}^2 E_*(\tilde{U}_{11}^{*i})^{\alpha+\beta}$$
$$+ 2\{E_* \tilde{Y}_{11}^{*i}(\tilde{U}_{11}^{*i})^\alpha E_* \tilde{Y}_{11}^{*i}(\tilde{U}_{11}^{*i})^\beta\}]$$
$$\equiv \frac{1}{\sqrt{n}} \hat{\chi}_{n,i}^{[3]} + \hat{Q}_{n,i}^{[3]}, \quad \text{say,}$$
$$\chi_*^{[4]}(T_{1n}^{*i}) = -\frac{2}{b\tilde{\tau}_{n,i}^4} E_*(\tilde{Y}_{11}^{*i})^4 + \frac{12}{b} + \hat{Q}_{n,i}^{[4]} \equiv b^{-1} \hat{\chi}_{n,i}^{[4]} + \hat{Q}_{n,i}^{[4]}, \quad \text{say,}$$

where, for some $N_5(\delta_1) \geq 1$ and $C(\delta_1) \in (0, \infty)$, on the set $A_{in}$, $\sum_{k=1}^{4} |\hat{Q}_{n,i}^{[k]} - [n^{-1}\Lambda_{4+k}(\{E_*(U_{21}^{*i})^\alpha : |\alpha| \leq 4\}; B_{n,i})]| \leq C(\delta_1) b^{-3/2}$ uniformly in $i = 0, 1, \ldots, M$ whenever $n \geq N_5(\delta_1)$. Hence the Edgeworth expansion for $T_{1n}^{*i}$ is given by

$$\hat{\Psi}_{n,i}(x) = \Phi(x) - \frac{1}{\sqrt{n}} \left( \hat{\chi}_{i,n}^{[1]} + \frac{1}{6} H_2(x) \hat{\chi}_{i,n}^{[3]} \right) \phi(x)$$
(5.14)
$$- b^{-1} \left( \frac{3}{2} H_1(x) + \frac{1}{24} H_3(x) \hat{\chi}_{n,i}^{[4]} \right) \phi(x)$$
$$+ n^{-1} \sum_{k=0}^{5} \Lambda_{9+K}(\{E_*(U_{21}^{*i})^\alpha : |\alpha| \leq 4\}; B_{n,i}) H_k(x) \phi(x),$$

where $H_k(x)$ denotes the $k$th order Hermite polynomial, $k \geq 1$.

We now complete the proof of the lemma by deriving the desired bound on $P(A_{in}^c)$, $0 \leq i \leq M$. By Lemmas 5.1 and 5.2, for all $n \geq N_6(\kappa, \delta_1)$,

$$P(A_{in}^{(1)c}) \leq P(\|\hat{\mu}_{n,i} - \mu\| > n^{-1/2} \log n)$$
$$+ \sum_{|\alpha|=3} P(\sqrt{\ell}|E_*(U_{11}^{*i})^\alpha - E_* U_{11}^\alpha| > c(\delta_1))$$
$$+ \sum_{|\alpha|=2}^{10} P(|E_*(U_{11}^{*i})^\alpha - E_* U_{11}^\alpha| > C(\delta_1))$$
$$\leq C(\delta_1, \rho(r, \delta))[o([n^{-1}m]^{5r}) + (\sqrt{\ell})^{6r}(n^{-1}\ell)^{3r} + n^{-r}\ell^r]$$
$$\leq C(\delta_1, \rho(r, \delta)) n^{-r} \ell^\gamma$$

uniformly in $i \in \{0, 1, \ldots, M\}$.

Also, by Theorem 2.8 of [15], for any $\delta_2 > 0, K > 0$, there exists $N_7(\kappa) \geq 1$ such that, for all $n \geq N_7(\kappa)$, $\sup\{|E \exp(it'U_{21})| : \delta_2 \leq \|t\| \leq \ell^K\} \leq \sup\{|E \exp(it'Z_{2,\infty})| : \delta_2 \leq \|t\| \leq \ell^K\} + C(\delta_2, K)\ell^{-1/2} \log \ell$ uniformly in $n^\kappa \leq \ell \leq \kappa^{-1} n^{1/3}$. Hence, by Lemma 5.2 (with $\lambda = \max\{8, r/2\}$) there exists an integer $N_8 = N_8(r, \kappa, \delta_1, \xi_\infty) \geq 1$ such that, for all $n \geq N_8$,

$$P(\hat{\xi}_{n,i} > \xi_\infty + 2^{-1}(1 - \xi_\infty), E_*\|U_{21}^*\| < C(\delta_1))$$
$$\leq P(\sup\{|E_* \exp(\sqrt{-1}t'U_{21}^*) - E \exp(\sqrt{-1}t'U_{21})| : \|t\| \leq C(\delta_1, k)n^8\}$$
$$> (1 - \xi_\infty)/4, E_*\|U_{21}^*\| < C(\delta_1))$$
$$\leq C(r, \theta_\infty)n^{-r}.$$

Therefore, it follows that $P(A_{in}^c) \leq C(r, \kappa, \delta_1, \xi_\infty)n^{-r}\ell^r$ for all $n \geq \max\{N_j : j = 1, \ldots, 8\}$, proving the lemma. □

PROOF OF THEOREM 3.1. First we simplify the leading terms in the Edgeworth expansion of $T_n^{*i}$ further. By Taylor's expansion of $\hat{c}_{\alpha,i} = D^\alpha H(\hat{\mu}_{n,i})/\alpha!$ around $\mu$, on the set $A_{in}$, we have

$$E_*(\tilde{Y}_{11}^{*i})^3 = E_*\left(\sum_{|\alpha|=1} \hat{c}_{\alpha,i}((U_{11}^{*i})^\alpha - E_*U_{11}^{*i})\right)^3$$

(5.15)
$$= E_*[\Gamma_{11}^{*i}]^3 - 3(E_*\Gamma_{11}^{*i})(E_*(\Gamma_{11}^{*i})^2)$$
$$+ \sum_{|\alpha|=2}\sum_{|\beta|=1} C(\{c_\gamma : |\gamma| \leq 2\})E_*(U_{11}^{*i})^{\alpha+\beta}\hat{\mu}_{n,i} + \hat{R}_{11n}^{(i)},$$

where $|\hat{R}_{11n}^{(i)}|\mathbb{1}_{A_{in}} \leq Cn^{-1}\sqrt{\ell}(\log n)^2$ for all $0 \leq i \leq M$. Next define the sets $A_{in}^{(3)} = \{|\tilde{\tau}_{n,i}^2 - \tau_\ell^2| \leq b^{-2/5}\}$ and $D_{i,n} = (A_{0n} \cap A_{in}) \cap (A_{0n}^{(3)} \cap A_{1n}^{(3)})$, $0 \leq i \leq M$. Then, using Lemma 5.1 and Chebyshev's inequality, one can show that, under the conditions of Lemma 5.4,

(5.16)
$$\max_{0 \leq i \leq M} P(D_{i,n}^c) \leq C(n^{-1}\ell)^r.$$

Next, using Taylor's expansion of $\tilde{\tau}_{i,n}^{-a} = (\tilde{\tau}_{i,n}^2)^{-a/2}$ around $\tau_\ell^2$ and (5.15), after some algebra one can show that the Edgeworth expansion of $T_n^{*i}$ may be rewritten as

$$\hat{\Psi}_{n,i}(x) = \Phi(x) + \frac{1}{\sqrt{n}}\frac{2x^2+1}{6\tau_\ell^3}\sqrt{\ell}[E_*(\Gamma_{11}^{*i})^3 - 3\tau_\ell^2 E_*\Gamma_{11}^{*i}]\phi(x)$$

$$+ \frac{1}{\sqrt{n}}\Lambda_{15}(\{\hat{U}_{1n}^{(i)}(\alpha) : 1 \leq |\alpha| \leq 3\}, \{\hat{\mu}_{n,i}\}, B; x)\phi(x)$$

(5.17)
$$+ b^{-1}\Lambda_{16}(\{\hat{U}_{1n}^{(i)}(\alpha) : 1 \leq |\alpha| \leq 4\}, B; x)\phi(x)$$

$$+ n^{-1}\Lambda_{17}(\{\hat{U}_{1n}^{(i)}(\alpha) : 1 \leq |\alpha| \leq 8\}, B; x)\phi(x) + \hat{R}_{12n}^{(i)}(x)$$

$$\equiv \hat{\Psi}_{1n}^{(i)}(x) + \hat{R}_{12n}^{(i)}(x), \quad \text{say,}$$

where $B = \{c_\alpha : |\alpha| \leq 3\} \cup \{\tau_\ell^{-1}\} \cup \{EU_{11}^\alpha : 1 \leq |\alpha| \leq 8\}$, and on the set $D_{i,n}$, $|\hat{R}_{12n}^{(i)}(x)| \leq C[b^{-3/2} + n^{-3/2}\ell](\log n)^4$ for all $0 \leq i \leq M, x \in \mathbb{R}$. Further, each term of the polynomial $\Lambda_{15}$ involves the product of at least two terms from the collection $\{\hat{U}_{1n}^{(i)}(\alpha) : 1 \leq |\alpha| \leq 3\} \cup \{\hat{\mu}_{n,i}\}$.

Next we find an expansion for $\text{Var}(\hat{\varphi}_{1n})$. Set $\Upsilon_{1n} = \Phi(x) + \frac{\sqrt{\ell}}{\sqrt{n}} \frac{(2x^2+1)}{6\tau_\ell^3} \times [E_*(\Gamma_{11}^{*0})^3 - 3\tau_\ell^2(E_*\Gamma_{11}^{*0})]\phi(x)$. Then, by Lemmas 5.1, 5.4 and (5.16),

$$
\begin{aligned}
(5.18) \quad &|\text{Var}(\hat{\varphi}_{1n}) - \text{Var}(\hat{\Upsilon}_{1n}\mathbb{1}_{D_{0,n}})| \\
&\leq 3P(D_{0,n}^c) + \text{Var}([\hat{\varphi}_{1n} - \hat{\Upsilon}_{1n}]\mathbb{1}_{D_{0,n}}) \\
&\quad + 2|\text{Cov}(\hat{\Upsilon}_{1n}\mathbb{1}_{D_{0,n}}, (\hat{\varphi}_{1n} - \hat{\Upsilon}_{1n})\mathbb{1}_{D_{0,n}})| \\
&\leq C[P(D_{0,n}^c) + \{\text{Var}(\hat{\Upsilon}_{1n}\mathbb{1}_{D_{0,n}})\}^{1/2}\{n^{-2}\ell\}^{1/2}].
\end{aligned}
$$

Using the arguments develped by Hall, Horowitz and Jing [20], one can show that

$$
(5.19) \quad \text{Var}(\hat{\Upsilon}_{1n}) = C\tau_\ell^{-6}(2x^2+1)^2\phi^2(x)n^{-2}\ell^2(1+o(1)).
$$

Also, by Hölder's inequality and Lemmas 5.1 and 5.4,

$$
\begin{aligned}
E(\hat{\Upsilon}_{1n}\mathbb{1}_{D_{0,n}^c})^2 &\leq 2(E(\hat{\Upsilon}_{1n} - E\hat{\Upsilon}_{1n})^2\mathbb{1}_{D_{0,n}^c} + (E\hat{\Upsilon}_{1n})^2 P(D_{0,n}^c)) \\
&\leq C[\{E(\hat{\Upsilon}_{1n} - E\hat{\Upsilon}_{1n})^4 P(D_{0,n}^c)\}^{1/2} + (n^{-1})P(D_{0,n}^c)] \\
&\leq C[\{n^{-1}\ell(n^{-1}\ell)^2(n^{-1}\ell)^3\}^{1/2} + n^{-1}n^{-3}\ell^3] \leq Cn^{-3}\ell^3.
\end{aligned}
$$

This implies that $|\text{Var}(\hat{\Upsilon}_{1n}\mathbb{1}_{D_{0,n}}) - \text{Var}(\hat{\Upsilon}_{1n})| \leq C(n^{-5}\ell^5)^{1/2}$ and, hence, from (5.18) and (5.19),

$$
(5.20) \quad |\text{Var}(\hat{\phi}_{1n}) - \text{Var}(\hat{\Upsilon}_{1n})| \leq Cn^{-2}\ell^{3/2}.
$$

Next we prove convergence of the JAB variance estimator. Let $W_{kj} = U_{1j}^{\alpha_k} - EU_{1j}^{\alpha_k}$ for some $\alpha_k \in \mathbb{Z}_+^d$ with $1 \leq |\alpha_k| \leq 8$, $k = 1, \ldots, r$, $1 \leq j \leq N$. Also, let $S_n(k,i) = |I_i|^{-1} \sum_{j \in I_i} W_{kj}$. Then, for $1 \leq i \leq M$,

$$
\begin{aligned}
\tilde{S}_{n,r}^{(i)} &\equiv \left\{m^{-1}\left[N\prod_{k=1}^r S_n(k,0) - (N-m)\prod_{k=1}^r S_n(k,i)\right]\right\} \\
&= \left[\sum_{k=1}^r \left\{\prod_{p=1}^{k-1} S_n(p,0)\right\}\left\{\prod_{p=k+1}^r S_n(p,i)\right\}\left\{m^{-1}\sum_{j \notin I_i} W_{kj}\right\}\right].
\end{aligned}
$$

Note that, by Lemma 5.1, $\max\{ES_n(k,i)^2 : 0 \leq i \leq M, 1 \leq k \leq r\} \leq Cn^{-1}\ell$, and by similar arguments, $\max\{E(m^{-1}\sum_{j \notin I_i} W_{kj})^2 : 1 \leq i \leq M, 1 \leq k \leq r\} \leq$ $Cm^{-1}\ell$. Further, on the set $D_{i,n}, |S_n(k,i)| \leq C$ for all $0 \leq i \leq M, 1 \leq k \leq r$. Hence, using the Cauchy–Schwarz inequality, for $r \geq 2$ we have

$$\left| M^{-1} \sum_{i=1}^{M} \tilde{S}_{n,r}^{(i)} \mathbb{1}_{D_{i,n}} \right|$$

$$\leq \sum_{k=1}^{r} \left\{ \prod_{p=1}^{k-1} S_n(p,0) \right\} \left( M^{-1} \sum_{i=1}^{M} \prod_{p=k+1}^{r} S_n^2(p,i) \mathbb{1}_{D_{i,n}} \right)^{1/2}$$

(5.21)
$$\times \left( M^{-1} \sum_{i=1}^{M} \left( m^{-1} \sum_{j \notin I_i} W_{kj} \right)^2 \right)^{1/2}$$

$$\leq \sum_{k=1}^{r} \left\{ \prod_{p=1}^{k-1} S_n(p,0) \right\}$$

$$\times \left( C(r) M^{-1} \sum_{i=1}^{M} S_n^2(p,i) \right)^{1/2} \left( M^{-1} \sum_{i=1}^{M} \left( m^{-1} \sum_{j \notin I_i} W_{kj} \right)^2 \right)^{1/2}$$

$$= O_p((n^{-1}\ell)^{1/2}(m^{-1}\ell)^{1/2}).$$

And for $r = 1$ there exist weights $\omega_{jn} \in [0,1], 1 \leq j \leq N$, such that

(5.22)
$$M^{-1} \sum_{i=1}^{M} \tilde{S}_{n,r}^{(i)} = M^{-1} \sum_{i=1}^{M} \left( m^{-1} \sum_{j \notin I_i} W_{kj} \right)$$

$$= M^{-1} \sum_{j=1}^{N} \omega_{jn} W_{kj} = O_p(n^{-1}\ell),$$

by arguments similar to Lemma 5.1. Let $\widehat{\text{VAR}}_{1n} = (m/N)M^{-1}\sum_{i=1}^{M}(\tilde{\Psi}_{1n}^{(i)}(x) - \hat{\Psi}_{1n}(x))^2 \mathbb{1}_{D_{i,n}}$ and $\hat{\Delta}_{1n} \equiv (m/N)M^{-1}\sum_{i=1}^{M}[(\tilde{\varphi}_{1n}^{(i)} - \hat{\varphi}_{1n}) - (\tilde{\Psi}_{1n}^{(i)}(x) - \hat{\Psi}_{1n}(x))]^2 \times \mathbb{1}_{D_{i,n}}$, where $\tilde{\Psi}_{1n}^{(i)}(x) = m^{-1}[N\hat{\Psi}_{1n}(x) - (N-m)\hat{\Psi}_{1n}^{(i)}(x)]$, $1 \leq i \leq M$. Then, using the fact that $\tilde{\varphi}_{1n}^{(i)} - \hat{\varphi}_{1n} = [(N-m)/m](\hat{\varphi}_{1n} - \hat{\varphi}_{1n}^{(i)})$, $1 \leq i \leq M$,

$$|\widehat{\text{Var}}_{\text{JAB}}(\hat{\varphi}_{1n}) - \widehat{\text{VAR}}_{1n}| \leq [\hat{\Delta}_{1n} + 2(\widehat{\text{VAR}}_{1n})^{1/2}(\hat{\Delta}_{1n})^{1/2}]$$

(5.23)
$$+ \frac{m}{N} M^{-1} \sum_{i=1}^{M} \left( \frac{N-m}{m} \right)^2 (\hat{\varphi}_{1n} - \hat{\varphi}_{1n}^{(i)})^2 \mathbb{1}_{D_{in}^c}.$$

By Lemma 5.4 and (5.17), it readily follows that

$$\hat{\Delta}_{1n} \leq \frac{m}{N} \left( \frac{N-m}{m} \right)^2 M^{-1} \sum_{i=1}^{M} [(\hat{\varphi}_{1n} - \hat{\Psi}_{1n}(x)) - (\hat{\varphi}_{1n}^{(i)} - \hat{\Psi}_{1n}^{(i)})]^2 \mathbb{1}_{D_{1n}}$$

$$
\begin{align}
(5.24) \quad &\leq C\frac{n}{m}[(n^{-3/2}\ell + b^{-3/2})(\log n)^4]^2 \\
&\leq Cn^{-2}\ell^2[m^{-1}\ell(\log n)^8].
\end{align}
$$

Since $|\hat{\varphi}_{1n}^{(i)} - \hat{\varphi}_{1n}| \leq 1$ for all $i$ and by Lemma 5.4, $M^{-1}\sum_{i=1}^m P(D_{in}^c) \leq C[n^{-1}\ell]^3$, by (5.18)–(5.20), (5.23) and (5.24), it is now enough to show that

$$
(5.25) \qquad |\widehat{\mathrm{VAR}}_{1n} - \mathrm{Var}(\hat{\Upsilon}_{1n})| = o_p(n^{-2}\ell^2).
$$

To that end, define the variables $\Gamma_{2j} = [n^{-1}\ell]^{1/2}(6\tau_i^3)^{-1}(2x^2+1)\phi(x)[\Gamma_{1j}^3 - 3\tau_\ell^2 \Gamma_{1j}]$, $1 \leq j \leq N$. Then from (5.17) we have $\hat{\Psi}_{1n}^{(i)}(x) = \Phi(x) + |I_i|^{-1} \times \sum_{j\in I_i}\Gamma_{2j} + \hat{R}_{13n}^{(i)}(x)$, $0 \leq i \leq M$, where $\hat{R}_{13n}^{(i)}(x)$ is the sum of the terms involving $\Lambda_{15}, \Lambda_{16}$ and $\Lambda_{17}$ in (5.17). Write $\tilde{R}_{13n}^{(i)}(x) = m^{-1}(N\hat{R}_{13n}^{(i)}(x) - (N-m) \times \hat{R}_{13n}^{(0)}(x))$, $\tilde{\Gamma}_{3i} = m^{-1}\sum_{j\notin I_i}(\Gamma_{2j} - E\Gamma_{21})$, $1 \leq i \leq M$, and $\bar{\Gamma}_{3n} = N^{-1}\sum_{j=1}^N(\Gamma_{2j} - E\Gamma_{21})$. Then it follows that

$$
\tilde{\Psi}_{1n}^{(i)}(x) - \hat{\Psi}_{1n}^{(0)}(x) = \tilde{\Gamma}_{3i} - \bar{\Gamma}_{3n} + \tilde{R}_{13n}^{(i)}(x) - \hat{R}_{13n}^{(0)}, \qquad 1 \leq i \leq M.
$$

Using this identity and noting that $\mathrm{Var}(\hat{\Upsilon}_{1n}) = \mathrm{Var}(\bar{\Gamma}_{3n})$, one has

$$
\begin{align}
&|\widehat{\mathrm{VAR}}_{1n} - \mathrm{Var}(\hat{\Upsilon}_{1n})| \\
(5.26) \quad &\leq \left|\widehat{\mathrm{VAR}}_{1n} - \frac{m}{N}\cdot\frac{1}{M}\sum_{i=1}^M \tilde{\Gamma}_{3i}^2\right| + \left|\frac{m}{N}\cdot\frac{1}{M}\sum_{i=1}^M (\tilde{\Gamma}_{3i}^2 - E\Gamma_{31}^2)\right| \\
&\quad + \left|\frac{m}{N}E\tilde{\Gamma}_{31}^2 - \mathrm{Var}(\bar{\Gamma}_{3n})\right| \\
&\equiv I_{11} + I_{12} + I_{13}, \qquad \text{say.}
\end{align}
$$

Clearly, $I_{13} = N^{-1}|\mathrm{Var}(m^{-1/2}\sum_{j=1}^m \Gamma_{2j}) - \mathrm{Var}(N^{-1/2}\sum_{j=1}^N \Gamma_{2j})| = o(n^{-2}\ell^2)$. Also, by arguments similar to the proof of Lemma 5.1, for any $0 < \eta < 1$,

$$
\begin{align}
(5.27) \quad EI_{12}^2 &\leq C \cdot \left[\frac{m}{N}\right]^2 M^{-2}\left[\frac{M}{\ell+m}\max_{1\leq k\leq \ell+m}\left(E\left|\sum_{i=1}^k \tilde{\Gamma}_{3i}^2\right|^{2+\eta}\right)^{2/(2+\eta)}\right] \\
&\leq C\cdot\left[\frac{m}{N}\right]^2 M^{-2}\left[\frac{M}{\ell+m}(\ell+m)^2(E|\tilde{\Gamma}_{31}|^{4+2\eta})^{2/(2+\eta)}\right] \\
&\leq Cm^3 n^{-3}\max_{|\alpha|=1,3}\left\{E\left|m^{-1}\sum_{j=1}^m\left(\frac{\ell}{n}\right)^{1/2}(U_{1j}^\alpha - EU_{11}^\alpha)\right|^{4+2\eta}\right\}^{2/(2+\eta)} \\
&\leq Cm^3 n^{-3}\left(\frac{\ell}{n}\right)^2 (m^{-1}\ell)^2 = Cn^{-5}m\ell^4.
\end{align}
$$

Next let $\hat{\Delta}_{2n} = \frac{m}{N}[m^{-1}\sum_{i=1}^{M} \tilde{R}_{13n}^2 \mathbb{1}_{D_{in}} + \bar{\Gamma}_{3n}^2 + \{\hat{R}_{13n}^{(0)}\}^2]$. Then it is easy to show that

$$I_{11} \leq C\left[\hat{\Delta}_{2n} + (\hat{\Delta}_{2n})^{1/2}\left\{\frac{m}{N}M^{-1}\sum_{i=1}^{M}\tilde{\Gamma}_{3i}^2\right\}^{1/2}\right]$$

(5.28)

$$+ C \cdot \frac{m}{N}M^{-1}\sum_{i=1}^{M}\tilde{\Gamma}_{3i}^2 \mathbb{1}_{D_{in}^c}.$$

By (5.19), (5.20), (5.26), (5.27) and the bound on $I_{13}$, $\frac{m}{N} \cdot M^{-1}\sum_{i=1}^{M}\tilde{\Gamma}_{3i}^2 = O_p(n^{-2}\ell^2)$. Also, using arguments similar to (5.27), one gets

$$E\left\{\frac{m}{N}M^{-1}\sum_{i=1}^{M}\tilde{\Gamma}_{3i}^2 \mathbb{1}_{D_{in}^c}\right\} \leq C\frac{m}{N} \cdot M^{-1}\sum_{i=1}^{M}[E\Gamma_{3i}^4]^{1/2}[P(D_{in}^c)]^{1/2}$$

(5.29)

$$\leq C\frac{m}{N}\left[\left(\frac{\ell}{n}\right)^2\left(\frac{\ell}{m}\right)^2\right]^{1/2}[n^{-3}\ell^3]^{1/2}$$

$$= o(n^{-2}\ell^2).$$

Finally, using (5.22), (5.23) and Lemma 5.1, one can show that $\hat{\Delta}_{2n} = o_p(n^{-2}\ell^2)$. Hence, (5.25) now follows from (5.26), (5.28) and (5.29). This completes the proof of Theorem 3.1. □

PROOF OF THEOREM 3.2. The proof is very similar to the proof of Theorem 3.1. Hence, we point out the necessary modifications. Using the arguments in the proof of Theorem 20.1 of [3], one can derive a version of Lemma 5.3 with a remainder term of the order $o(n^{-3/2})$. Next, using arguments similar to the proof of Lemma 5.4, on the set $A_{in}$, one can derive a fourth-order Edgeworth expansion $\hat{\Psi}_{4n}^{(i)}$, say, for the Studentized statistic $T_n^{*i}$ with an error term $o(n^{-3/2})$ uniformly in $0 \leq i \leq n$. Here the $O(n^{-1/2})$ terms of $\hat{\Psi}_{4n}^{(i)}$ are the same as those given by Lemma 5.4 with $\ell = 1 = m$. Hence, it follows that $\text{Var}(\hat{\varphi}_{1n}) \sim Cn^{-2}$ for all $x_0 \in \mathbb{R}$, which goes to zero at a rate faster than $\text{Var}(\hat{\varphi}_{1n})$ in the dependent case. This is why we needed a more accurate Edgeworth expansion in the independent case than the dependent case.

Next note that, by Corollary 4.4 of [12], for any i.i.d. random variables $W_1, \ldots, W_n$ with $EW_1 = 0$ and $E|W_1|^t < \infty$ for some $t \in [2, \infty)$,

$$P\left(\sum_{i=1}^{n}|W_i| \geq x\right) \leq C_1(t)[nE|W_1|^t]x^{-t}$$

(5.30)

$$+ \exp(-C_2(t)x^2/[n\text{Var}(W_1)]), \qquad x > 0.$$

Now using (5.30) and moderate deviation probability inequalities (cf. [3]) for sums of independent random variables in place of Lemmas 5.1 and 5.2, one can complete the proof of Theorem 3.2 by retracing the steps in the proof of Theorem 3.1. □

PROOF OF THEOREM 3.3. Let $z_\alpha$ denote the $\alpha$-quantile of $N(0,1)$. Then, inverting the Edgeworth expansion of $T_n^{*i}$, given by Lemma 5.4 and (5.17), we have

$$
\begin{aligned}
\hat{\varphi}_{2n}^{(i)} = {} & z_\alpha - \frac{\sqrt{\ell}}{\sqrt{n}} \frac{(2z_\alpha^2+1)}{6\tau_\ell^3}[E_*(\Gamma_{11}^{*i})^3 - 3\tau_\ell^2 E_*\Gamma_{11}^{*i}] \\
& + n^{-1/2}\Lambda_{18}(\{\hat{U}_{1n}^{(i)}(\alpha): 1 \leq |\alpha| \leq 3\}, \{\hat{\mu}_{n,i}\}, B; z_\alpha) \\
& + b^{-1}\Lambda_{19}(\{\hat{U}_{1n}^{(i)}(\alpha): 1 \leq |\alpha| \leq 4\}, z_\alpha) \\
& + n^{-1/2}b^{-1/2}\Lambda_{20}(\{\hat{U}_{1n}^{(i)}(\alpha): 1 \leq |\alpha| \leq 4\}, \{\hat{\mu}_{n,i}\}, B; z_\alpha) \\
& + n^{-1}\Lambda_{21}(\{\hat{U}_{1n}^{(i)}(\alpha): 1 \leq |\alpha| \leq 8\}, B; z_\alpha) + \hat{R}_{21n}^{(i)} \\
\equiv {} & \hat{\Psi}_{2n}^{(i)} + \hat{R}_{21n}^{(i)}, \qquad \text{say,}
\end{aligned}
$$
(5.31)

where, on the set $D_{in}, |\hat{R}_{21n}^{(i)}| \leq C[n^{-3/2}\ell + b^{-3/2}](\log n)^4$ for all $0 \leq i \leq M$. Write $\hat{\Upsilon}_{2n}^{(i)} = z_\alpha - \frac{\sqrt{\ell}}{\sqrt{n}}\frac{2z_\alpha^2+1}{6\tau_\ell^3}[E_*(\Gamma_{11}^{*i})^3 - 3\tau_\ell^2 E_*\Gamma_{11}^{*i}], 0 \leq i \leq M$. Then

$$
\begin{aligned}
|\operatorname{Var}(\hat{\varphi}_{2n}) - \operatorname{Var}(\hat{\Upsilon}_{2n})| \leq {} & E(\hat{\varphi}_{2n} - \hat{\Upsilon}_{2n})^2 \\
& + 2[E(\hat{\varphi}_{2n} - \hat{\Upsilon}_{2n})^2]^{1/2}[\operatorname{Var}(\hat{\Upsilon}_{2n})]^{1/2}.
\end{aligned}
$$
(5.32)

Note that for any random variable $W$ and $0 < \alpha < 1$, if $P(|W| > C_0) < \min\{\alpha, 1-\alpha\}$ for some $C_0 \in (0, \infty)$, then its $\alpha$-quantile $w_\alpha$ satisfies the inequality $|w_\alpha| \leq C_0$. Since $\|\hat{\mu}_{n,i}\| \leq \ell^{-1/2}E_*\|U_{11}^{*i}\|$, by Taylor's expansion and assumption (A.1), for $0 \leq i \leq M$ and $t_0 \in (0, \infty)$ we have

$$
\begin{aligned}
P_*(|T_n^{*i}| > t_0) & \leq t_0^{-2}E_*(T_n^{*i})^2 \leq t_0^{-2}n^3 E_*[H(\bar{X}_n^{*i}) - H(\hat{\mu}_{n,i})]^2 \\
& \leq t_0^{-2}n^3 C(d,a)\bigg[ E_*\bigg\{\sum_{|\alpha|=1}\hat{c}_{\alpha,i}(\bar{X}_n^{*i} - \hat{\mu}_{n,i})^\alpha\bigg\}^2 \\
& \qquad\qquad + E_*\{(1 + \|\bar{X}_n^{*i}\|^a + \|\hat{\mu}_{n,i}\|^a)\|\bar{X}_n^{*i} - \hat{\mu}_{n,i}\|^2\}^2\bigg] \\
& \leq t_0^{-2}n^3 C(d,\kappa,a)[n^{-1}\tilde{\tau}_{n,i}^2 + n^{-2}\{1 + E_*\|U_{11}^{*i}\|^{4+2a}\}] \\
& \leq t_0^{-2}n^2 C(d,\kappa,a)[1 + E_*\|U_{11}^{*i}\|^{4+2a}] < \min\{\alpha, 1-\alpha\},
\end{aligned}
$$

provided $t_0 \geq Cn\{1 + E_*\|U_{11}^{*i}\|^{4+2a}\}^{1/2}$ for some $C = C(d, \kappa, a, \alpha) \in (0, \infty)$ independent of $i$. Hence it follows that

(5.33) $\quad \hat{\varphi}_{2n}^{(i)} \leq Cn\{1 + E_*\|U_{11}^{*i}\|^{4+2a}\}^{1/2} \quad$ for all $0 \leq i \leq M$.

Now, using Hölder's inequality and Lemmas 5.1 and 5.4, we get

$$
\begin{aligned}
&E(\hat{\varphi}_{2n} - \hat{\Upsilon}_{2n})^2 \\
&\leq E\hat{\varphi}_{2n}^2 \mathbb{1}_{D_{0,n}^c} + E\Upsilon_{2n}^2 \mathbb{1}_{D_{0,n}^c} + E(\hat{\varphi}_{2n} - \hat{\Upsilon}_{2n})^2 \mathbb{1}_{D_{0,n}} \\
&\leq Cn^2[(1 + E\|U_{11}\|^{4+2a})P(D_{0,n}^c) \\
&\qquad + \{E(E_*\|U_{11}^*\|^{4+2a} - E\|U_{11}\|^{4+2a})^{12}\}^{1/12}\{P(D_{0,n}^c)\}^{11/12}] \\
&\quad + C(D_H, z_\alpha, \kappa)b^{-1} \max\{E|E_*(\Gamma_{11}^{*0})^k|\mathbb{1}_{D_{0,n}^c} : k = 1, 3\} \\
&\quad + C[n^{-3}\ell^2 + b^{-3}](\log n)^8 \\
&\leq Cn^2[P(D_{0,n}^c) + \{(n^{-1}\ell)^6\}^{1/12}\{P(D_{0,n}^c)\}^{11/12}] \\
&\quad + C(D_H, z_\alpha, \kappa)b^{-1} \max\{(E\Gamma_{11}^{2k})^{1/2}[P(D_{0,n}^c)]^{1/2} : k = 1, 3\} \\
&\quad + C[n^{-3}\ell^2 + b^{-3}](\log n)^8 \\
&= o(n^{-2}\ell^2),
\end{aligned}
$$
(5.34)

which, in view of (5.32), implies that $\text{Var}(\hat{\varphi}_{2n}) \sim Cn^{-1}\ell^2(1 + o(1))$. Thus it is now enough to show that

(5.35) $\quad |\widehat{\text{Var}}_{\text{JAB}}(\hat{\varphi}_{2n}) - \text{Var}(\hat{\Upsilon}_{2n})| = o_p(n^{-2}\ell^2)$.

For this define the variables $\widehat{\text{VAR}}_{2n}$ and $\hat{\Delta}_{3n}$ by replacing $\hat{\Psi}_{1n}^{(i)}(x)$'s and $\hat{\varphi}_{1n}^{(i)}$'s in the definitions of $\widehat{\text{VAR}}_{1n}$ and $\hat{\Delta}_{1n}$ by $\hat{\Psi}_{2n}^{(i)}$'s and $\hat{\varphi}_{2n}^{(i)}$'s, respectively. Then, by (5.31), $\hat{\Delta}_{3n} \leq Cn^{-2}\ell^2(m^{-1}\ell(\log n)^8)$, as in (5.24). Also, writing $U_{3i}(a) \equiv E_*\|U_{11}^{*i}\|^{4+2a} - E\|U_{11}\|^{4+2a}$, we have $E|U_{3i}(a)|\mathbb{1}_{D_{i,n}^c} \leq (E|U_{3i}(a)|^{12})^{1/12} \times [P(D_{i,n}^c)]^{11/12} \leq Cn^{-6}\ell^6$ uniformly in $1 \leq i \leq M$. Hence, by (5.16) and (5.33),

$$
\frac{m}{N}M^{-1}\sum_{i=1}^{M}\left(\frac{N-m}{m}\right)^2(\hat{\phi}_{2n} - \hat{\phi}_{2n}^{(i)})^2 \mathbb{1}_{D_{in}^c}
$$

$$
\leq Cm^{-1}n^3M^{-1}\sum_{i=1}^{M}\{(1 + E_*\|U_{11}^*\|^{4+2a} + E_*\|U_{11}^{*i}\|^{4+2a}\}\mathbb{1}_{D_{i,n}^c}
$$

(5.36) $\quad \leq C\dfrac{n^3}{m}\left[(1 + E_*\|U_{11}^*\|^{4+2a} + E\|U_{11}\|^{4+2a})\left\{\dfrac{1}{M}\sum_{i=1}^{M}\mathbb{1}_{D_{i,n}^c}\right\}\right.$

$$+ \frac{1}{M}\sum_{i=1}^{M} |U_{3i}(a)|\mathbb{1}_{D_{i,n}^c}\Bigg]$$

$$= Cm^{-1}n^3[O_p(1)O_p((n^{-1}\ell)^6) + O_p((n^{-1}\ell)^6)] = o_p(n^{-2}\ell^2).$$

Therefore, using (5.36), the bound on $\hat{\Delta}_{3n}$ and arguments similar to (5.23), one gets $|\widehat{\mathrm{Var}}_{\mathrm{JAB}}(\hat{\varphi}_{2n}) - \widehat{\mathrm{VAR}}_{2n}| = o_p(n^{-2}\ell^2)$. Now, retracing the steps in the proof of (5.25), one can show that $|\widehat{\mathrm{VAR}}_{2n} - \mathrm{Var}(\hat{\Upsilon}_{2n})| = o_p(n^{-2}\ell^2)$, which yields (5.35). This completes the proof of Theorem 3.3. $\square$

PROOF OF THEOREM 3.4. Note that in the i.i.d. case, $\ell = m \equiv 1$ implies that $N = M = n$ and hence, $U_{3i}(a) \equiv E_*\|U_{11}^{*i}\|^{4+2a} - E\|U_{11}\|^{4+2a} = |I_i|^{-1}\sum_{j \in I_i}(\|X_j\|^{4+2a} - E\|X_1\|^{4+2a})$, where $I_0 = \{1,\ldots,n\}$ and $I_i = I_0 \setminus \{i\}$ for $1 \leq i \leq M$. Since the $X_i$'s are i.i.d., $U_{3i}(a)$ has the same distribution as $U_{31}(a)$ for all $1 \leq i \leq M$. Now, using (5.30) with $t = 7.5$ and using Theorem 3.21 of [14], we have

$$\begin{aligned}
E|U_{3i}(a)|\mathbb{1}_{D_{i,n}^c} &\leq \{E|U_{3i}(a)|^4\mathbb{1}(|U_{3i}(a)| > |I_i|^{-1/2}\log n)\}^{1/4}[P(D_{in}^c)]^{3/4} \\
&\quad + \{|I_i|^{-1/2}\log n\}P(D_{i,n}^c) \\
&\leq \{o(n^{-1})\}^{1/4}\{Cn^{-6.5}\}^{3/4} + C\{n^{-1/2}\log n\}n^{-6.5} \\
&= o(n^{-5}),
\end{aligned}$$

uniformly in $0 \leq i \leq M$. The proof of the theorem can now be completed by combining the arguments in the proofs of Theorems 3.2 and 3.3 and using the above inequality in verification of (5.34) and (5.36). We omit the routine details.

$\square$

**Acknowledgments.** The author thanks an Associate Editor and two referees for constructive suggestions that improved an earlier draft of the paper. Special thanks are due to Mr. K. Furukawa for his help with the simulation results.

DEPARTMENT OF STATISTICS
IOWA STATE UNIVERSITY
AMES, IOWA 50011-3440
USA
E-MAIL: snlahiri@iastate.edu